\documentclass[12pt]{amsart}

\usepackage{latexsym, amsmath, amssymb, amsfonts, amscd, amsthm, mathrsfs}

\newcommand{\Ce}{{\mathbb C}}
\renewcommand{\Re}{{\mathbb R}}

\newcommand{\Ze}{{\mathbb Z}}
\newcommand{\Ne}{{\mathbb N}}
\newcommand{\Te}{{\mathbb T}}

\newcommand{\cT}{{\mathcal T}}
\newcommand{\cA}{{\mathcal A}}
\newcommand{\cB}{{\mathcal B}}

\newcommand{\Der}{{\rm Der}}
\newcommand{\ADer}{{\rm ADer}}
\newcommand{\Hom}{{\rm Hom}}
\newcommand{\End}{{\rm End}}
\newcommand{\SO}{{\rm SO}}
\newcommand{\rO}{{\rm O}}
\newcommand{\GL}{{\rm GL}}
\newcommand{\ch}{{\rm ch}}
\newcommand{\id}{{\rm id}}
\newcommand{\Out}{{\rm Out}}
\newcommand{\rH}{{\rm H}}

\newcommand{\pa}{\|}

\newcommand{\ept}{\emptyset}

\theoremstyle{plain}
\newtheorem{theorem}{Theorem}[section]
\newtheorem{lemma}[theorem]{Lemma}
\newtheorem{proposition}[theorem]{Proposition}
\newtheorem{corollary}[theorem]{Corollary}
\newtheorem{question}[theorem]{Question}
\theoremstyle{definition}
\newtheorem{remark}[theorem]{Remark}
\newtheorem{example}[theorem]{Example}
\newtheorem{notation}[theorem]{Notation}
\newtheorem{definition}[theorem]{Definition}

\begin{document}

\title{Morita Equivalence of Smooth Noncommutative Tori}
\author{George A. Elliott}
\address{Department of Mathematics \\
University of Toronto \\
Toronto, Ontario, Canada\,\, M5S 2E4}
\email{elliott@math.toronto.edu}

\author{Hanfeng Li}
\address{Department of Mathematics \\
SUNY at Buffalo \\
Buffalo, NY 14260, USA} \email{hfli@math.buffalo.edu}

\date{July 7, 2006}


\begin{abstract}
We show that in the generic case the smooth noncommutative tori
associated to two $n\times n$ real skew-symmetric matrices are
Morita equivalent if and only if the matrices are in the same orbit
of the natural $\SO(n, n|\, \Ze)$ action.
\end{abstract}

\thanks{This research was supported by a grant from the
Natural Sciences and Engineering Research Council of Canada, held by
the first named author.}

\maketitle

\section{Introduction}
\label{intro:sec}

Let $n\ge 2$ and denote by $\mathcal{T}_n$ the space of $n\times n$ real
skew-symmetric matrices. For each $\theta \in \mathcal{T}_n$ the
corresponding $n$-dimensional noncommutative torus $A_{\theta}$ is
defined as the universal C*-algebra generated by unitaries
$U_1, {\cdots}, U_n$ satisfying the relation
\begin{eqnarray*}
U_kU_j=e(\theta_{kj})U_jU_k,
\end{eqnarray*}
where $e(t)=e^{2\pi it}$. Noncommutative tori are one of the
canonical examples in noncommutative differential geometry
\cite{Rieffel90, Connes94}.

One may also consider the smooth version $A^{\infty}_{\theta}$ of
a noncommutative torus, which is the algebra of formal series
\begin{eqnarray*}
\sum c_{j_1, {\cdots}, j_n}U^{j_1}_1\cdots U^{j_n}_n,
\end{eqnarray*}
where the coefficient function  $\Ze^n\ni(j_1, {\cdots}, j_n)\mapsto
c_{j_1, {\cdots}, j_n}$ belongs to the Schwartz space
$\mathcal{S}(\Ze^n)$. This is the space of smooth elements of
$A_{\theta}$ for the canonical action of $\mathbb{T}^n$ on
$A_{\theta}$.

A notion of Morita equivalence of C*-algebras (as an analogue of
Morita equivalence of unital rings \cite[Chapter 6]{AF74}) was
introduced by Rieffel in \cite{Rieffel74, Rieffel82}. This is now
often called Rieffel-Morita equivalence. It is known that two
unital C*-algebras are Morita equivalent as unital
$\Ce$-algebras if and only if
they are Rieffel-Morita equivalent \cite[Theorem 1.8]{Beer82}.
Rieffel-Morita equivalent C*-algebras share a lot in common such
as equivalent categories of Hilbert C*-modules, isomorphic {\rm
K}-groups, etc., and hence are usually thought of as having the
same geometry.

In \cite{Schwarz98} Schwarz introduced the notion
of complete Morita equivalence of smooth noncommutative tori,
which includes Rieffel-Morita equivalence of the corresponding
C*-algebras, but is stronger, and has
important application in M(atrix) theory \cite{Schwarz98, KS02}.

A natural question is to classify
noncommutative tori and their smooth counterparts up to
the various notions of Morita equivalence. Such results have
important application to physics  \cite{CDS98, Schwarz98}. For
$n=2$ this was done by Rieffel \cite{Rieffel81}. (For the
earlier problem of isomorphism, see below.)
In this case it does not
matter what kind of Morita equivalence we are referring to: there
is a (densely defined) action of the group $\GL(2, \Ze)$ on
$\mathcal{T}_2$, and two matrices in $\mathcal{T}_2$ give rise to
Morita equivalent noncommutative tori or smooth noncommutative
tori if and only if they are in the same orbit of this action, and
also if and only if the ordered ${\rm K}_0$-groups of the algebras are
isomorphic. The higher dimensional case is much more complicated
and there are examples showing that the smooth counterparts of two
Rieffel-Morita equivalent noncommutative tori may fail to be
completely Morita equivalent \cite{Rieffel99a}.

In \cite{Rieffel99a}
Rieffel and Schwarz found a (densely defined) action of the group
$\SO(n, n|\, \Ze)$ on $\mathcal{T}_n$ generalizing the above $\GL(2,
\Ze)$ action. Recall that  $\rO(n, n|\, \Re)$ denotes the group of
linear transformations of the vector space $\Re^{2n}$ preserving
the quadratic form $x_1x_{n+1}+x_2x_{n+2}+\cdots+x_nx_{2n}$, and
that $\SO(n, n|\, \Ze)$ refers to the subgroup of $\rO(n, n|\, \Re)$
consisting of matrices with integer entries and determinant $1$.

Following \cite{Rieffel99a}, let us write the elements of $\rO(n, n|\, \Re)$
in $2\times 2$ block form:
\begin{eqnarray*}
g=\begin{pmatrix} A & B \\ C &  D \end{pmatrix}.
\end{eqnarray*}
Then $A, B, C$, and $D$ are arbitrary $n\times n$ matrices satisfying
\begin{eqnarray} \label{O(n,n|R):eq}
A^tC+C^tA=0=B^tD+D^tB, & A^tD+C^tB=I.
\end{eqnarray}
The action of $\SO(n, n|\, \Ze)$ is then defined as
\begin{eqnarray} \label{action:eq}
g\theta=(A\theta+B)(C\theta+D)^{-1},
\end{eqnarray}
whenever $C\theta+D$ is invertible. There is a dense subset of
$\mathcal{T}_n$ on which the action of every $g\in \SO(n, n|\, \Ze)$
is defined \cite[page 291]{Rieffel99a}.

After the work of Rieffel, Schwarz, and the second named author in
\cite{Rieffel99a, Schwarz98, LiMorita} (see also \cite{TW}) it is
now known that two matrices in $\mathcal{T}_n$ give completely
Morita equivalent smooth noncommutative tori (in the sense of
\cite{Schwarz98}) if and only if they are in the same orbit of the
$\SO(n, n|\, \Ze)$ action.

Phillips has been able to
show that two simple noncommutative tori are Rieffel-Morita
equivalent if and only if their ordered ${\rm K}_0$-groups are
isomorphic \cite[Remark 7.9]{Phillips03}. Using the result in
\cite{LiMorita} and Phillips's result,
recently we have completed the classification of noncommutative
tori up to Rieffel-Morita equivalence \cite{EL2}. In general,
two noncommutative tori
are Rieffel-Morita equivalent if and only if they have isomorphic
ordered ${\rm K}_0$-groups and centers.

It remains to
classify smooth noncommutative tori up to Morita
equivalence as unital algebras. We shall consider a natural subset $\mathcal{T}'_n\subseteq
\mathcal{T}_n$ which can be described both algebraically
in terms of the properties of the algebra $A^{\infty}_{\theta}$
(see Notation~\ref{generic:notation} and Corollary~\ref{generic:corollary})
and number theoretically (see Proposition~\ref{generic 2:prop}),
and has the property that the complement
$\mathcal{T}_n\setminus \mathcal{T}'_n$ has Lebesgue measure zero
(Proposition~\ref{generic:prop}). The main result of this paper
is the Morita equivalence classification of the algebras arising from
the subset $\mathcal{T}'_n$:

\begin{theorem} \label{Morita equiv:thm}
(1) The set $\mathcal{T}'_n$ is closed under Morita equivalence of the associated
smooth noncommutative tori; in other words,
if the algebras $A^{\infty}_{\theta}$ and $A^{\infty}_{\theta'}$ are Morita equivalent and
$\theta\in \mathcal{T}'_n$, then $\theta'\in \mathcal{T}'_n$.

(2) Two matrices in $\mathcal{T}'_n$ give rise to Morita equivalent
smooth noncommutative tori if and only if they are in the same
orbit of the $\SO(n, n|\, \Ze)$ action.
\end{theorem}

Denote by $\mathcal{T}^{\flat}_n$ the subset of $\cT_n$ consisting
of $\theta$'s such that $A_{\theta}$ is simple. The weaker form of
the part (1) of Theorem~\ref{Morita equiv:thm} with $\mathcal{T}'_n$
replaced by $\mathcal{T}'_n\cap \mathcal{T}^{\flat}_n$ is a
consequence of \cite{Nest88} (see the discussion after the proof of
Proposition~\ref{generic 2:prop}).

 Consider the subset
$\tilde{\mathcal{T}}_3$ of $\mathcal{T}_3$ consisting of $\theta$'s
such that the seven numbers consisting of $1, \, \theta_{12}, \,
\theta_{13},$ $\theta_{23}$, together with all products of any two
of these four, are linearly independent over the rational numbers.
In \cite{Rieffel99a} Rieffel and Schwarz showed that for any
$\theta\in \tilde{\mathcal{T}}_3$, the matrices $\theta$ and
$-\theta$ are not in the same orbit of the $\SO(n, n|\, \Ze)$
action, although, by the work of Q. Lin and the first named author
on the structure of $3$-dimensional simple noncommutative tori
culminating in \cite{Lin96}, the C*-algebras $A_{\theta}$ and
$A_{-\theta}$ are isomorphic. It is easy to see that the complement
$\mathcal{T}_3\setminus \tilde{\mathcal{T}}_3$ has Lebesgue measure
zero. Our Theorem~\ref{Morita equiv:thm} (together with
Proposition~\ref{generic:prop}) shows that the complement
$\mathcal{T}_3\setminus (\mathcal{T}'_3\cap \tilde{\mathcal{T}}_3)$
also has Lebesgue measure zero, and for any $\theta\in
\mathcal{T}'_3\cap \tilde{\mathcal{T}}_3$, the algebras
$A^{\infty}_{\theta}$ and $A^{\infty}_{-\theta}$ are not Morita
equivalent.

A related question is the classification of noncommutative tori and
their smooth counterparts up to isomorphism. The case $n=2$ was done
by Pimsner, Rieffel and Voiculescu \cite{PV80, Rieffel81} and the
simple C*-algebra case for $n> 2$ was also done by Phillips
\cite[Theorem 7.6]{Phillips03} (see \cite[Section 5]{Phillips03} for
more of the history). There have been various results for the smooth
algebra case with $n>2$ \cite{DEKR85, CEGJ85, BCEN87}. In
particular, Cuntz, Goodman, Jorgensen, and the first named author
showed that two matrices in $\mathcal{T}'_n\cap
\mathcal{T}^{\flat}_n$ give isomorphic smooth noncommutative tori if
and only if the associated skew-symmetric bicharacters of $\Ze^n$
are isomorphic \cite{CEGJ85}.
This result is essentially a special case of
Theorem~\ref{Morita equiv:thm}, and the proof of it obtained in this
way (see Remark~\ref{isomorphism:remark}) is new.

Schwarz proved the ``only if'' part of Theorem~\ref{Morita
equiv:thm}(2) in the context of complete Morita equivalence
\cite[Section 5]{Schwarz98}. His proof is based on the Chern
character \cite{Connes80, Elliott84}, which is essentially a
topological algebra invariant. In order to show that his argument
still works in our situation, we have to show that a purely algebraic
Morita equivalence between smooth noncommutative tori is
automatically ``topological'' in a suitable sense.
For this purpose and also for the proof
of Theorem~\ref{Morita equiv:thm}(1), in Sections~\ref{cont of
alg:sec} and \ref{cont of der:sec} we show that any algebraic
isomorphism between two ``smooth algebras'' (see Remark~\ref{unique topology:remark} below)
is continuous
and any derivation of a ``smooth algebra'' is continuous. These
are the noncommutative analogues of the following well-known facts in classical
differential geometry: any algebraic isomorphism between the
smooth function algebras of two smooth manifolds corresponds to a
diffeomorphism between the manifolds, and any derivation of
the smooth function algebra corresponds to a (complexified) smooth
vector field on the manifold. We introduce the set
$\mathcal{T}'_n$ and prove Theorem~\ref{Morita equiv:thm}(1) in
Section~\ref{gs:sec}. Theorem~\ref{Morita equiv:thm}(2) will be
proved in Section~\ref{proof:sec}.

Throughout this paper $A$ will be a C*-algebra, and
$A^{\infty}$ will be a dense sub-$*$-algebra of $A$ closed under
the holomorphic functional calculus (after the adjunction of a unit) and
equipped with a
Fr\'echet space topology stronger than the C*-algebra norm topology.
Unless otherwise specified, the topology considered on $A^{\infty}$ will
always be this Fr\'echet topology.

We thank the referees and Ryszard Nest for very helpful comments.
This work was carried out while H. Li was at the University of
Toronto.

\section{Continuity of Algebra Isomorphisms}
\label{cont of alg:sec}

In this section we shall prove Proposition~\ref{topological alg:prop} and Theorem~\ref{cont of isom:thm}, which
indicate that the topology of $A^{\infty}$ necessarily behaves well with respect to the algebra structure.

\begin{lemma} \label{cont:lemma}
Let $\varphi:A^{\infty}\rightarrow A^{\infty}$ be an $\Re$-linear map.
If $A^{\infty}\overset{\varphi}\rightarrow A$ is continuous, then so
is $A^{\infty}\overset{\varphi}\rightarrow A^{\infty}$.
\end{lemma}
\begin{proof}
We shall use the closed graph theorem \cite[Corollary 48.6]{Berberian74}
to prove the continuity of $A^{\infty}\overset{\varphi}\rightarrow
A^{\infty}$. Since $A^{\infty}\overset{\varphi}\rightarrow A$ is continuous,
the graph of $\varphi$ is closed with respect to the norm topology in the second coordinate.
It is then also closed with respect to the Fr\'echet topology in the second coordinate.
It follows that $A^{\infty}\overset{\varphi}\rightarrow A^{\infty}$ is continuous.
\end{proof}

\begin{proposition} \label{topological alg:prop}
In $A^{\infty}$ the $*$-operation is continuous, and the multiplication is jointly continuous.
In other words, $A^{\infty}$ is a topological $*$-algebra.
\end{proposition}
\begin{proof}
By Lemma~\ref{cont:lemma}, the $*$-operation is continuous and
the multiplication is separately continuous in $A^{\infty}$.
Proposition~\ref{topological alg:prop} follows because of the fact
that if the multiplication of an algebra equipped with a Fr\'echet topology is separately continuous
then it is jointly continuous \cite[Proposition VII.1]{Waelbroeck71}.
\end{proof}

The following theorem was proved by Gardner in the case
$A^{\infty}=A$ \cite[Proposition 4.1]{Gardner65}, and was
given as Lemma 4 of \cite{CEGJ85} in the case of smooth
noncommutative tori.

\begin{theorem} \label{cont of isom:thm}
Every algebra isomorphism $\varphi:A^{\infty}_1\rightarrow A^{\infty}_2$ is
an isomorphism of topological algebras.
\end{theorem}
\begin{proof}
Our proof is a modification of that for \cite[Proposition 4.1]{Gardner65}.
It suffices to show that $\varphi^{-1}$ is continuous.
For $a\in A^{\infty}_1$ let $r(a)$ denote the spectral radius of $a$,
in $A^{\infty}_1$ and also in $A_1$ \cite[Lemma 1.2]{Schweitzer92}.

First, for any $a\in A^{\infty}_1$ we have
$$
\pa a\pa^2=\pa aa^*\pa=r(aa^*)=r(\varphi(a)\varphi(a^*))\le \pa
\varphi(a)\varphi(a^*)\pa \le \pa \varphi(a)\pa \cdot \pa \varphi(a^*)\pa.
$$

Next, we use the closed graph theorem (cf.~above) to show that
$\varphi\circ *\circ \varphi^{-1}$ is continuous on $A^{\infty}_2$.
Let $\{a_m\}_{m\in\Ne}\subseteq A^{\infty}_1$ be such that
$\varphi(a_m)\rightarrow 0$ and $\varphi(a^*_m)\rightarrow \varphi(b)$ for
some $b\in A^{\infty}_1$. By the preceding inequality we have
\begin{eqnarray*}
\pa a_m\pa^2\le \pa \varphi(a_m)\pa \cdot \pa \varphi(a^*_m)\pa \to 0 \cdot \pa \varphi(b)\pa=0, \\
\pa a^*_m-b\pa^2\le \pa \varphi(a^*_m)-\varphi(b)\pa \cdot \pa \varphi(a_m)-\varphi(b^*)\pa \to 0\cdot \pa \varphi(b^*)\pa=0.
\end{eqnarray*}
Therefore $b=0$. This shows that the graph of $\varphi\circ *\circ \varphi^{-1}$ is closed.

Finally, let us use the closed graph theorem to show that $\varphi^{-1}$ is continuous.
Let $\{a_m\}_{m\in\Ne}\subseteq A^{\infty}_1$ be such that
$\varphi(a_m)\rightarrow 0$ and $a_m\rightarrow b$ for
some $b\in A^{\infty}_1$. By continuity of $\varphi\circ *\circ \varphi^{-1}$ we have $\varphi(a^*_m)\to 0$.
By the inequality derived in the second paragraph it follows that $\pa a_m\pa^2\to 0$.
Therefore $b=0$. This shows that the graph of
$\varphi^{-1}$ is closed.
\end{proof}

\begin{question} \label{cont of isom:question}
Is every algebra isomorphism $\varphi:A^{\infty}_1\rightarrow A^{\infty}_2$ continuous with respect to the C*-algebra norm topology?
\end{question}

\begin{remark} \label{unique topology:remark}
Let us say that a Fr\'echet topology on an algebra $\mathcal{A}$ is a {\it smooth topology}
if there is a C*-algebra $A$ and a continuous injective homomorphism $\varphi:\mathcal{A}\hookrightarrow
A$ such that $\varphi(\mathcal{A})$ is a dense sub-$*$-algebra of $A$ closed under the holomorphic functional calculus.
Theorem~\ref{cont of isom:thm} can be restated as that {\it every algebra admits at most one smooth topology}.
\end{remark}

\begin{example} \label{cut down:example}
Let $p\in M_n(A^{\infty})$ be a projection. Then $pM_n(A^{\infty})p$ is a
dense sub-$*$-algebra of $pM_n(A)p$, and the relative topology on $pM_n(A^{\infty})p$
is a Fr\'echet topology stronger than the C*-algebra norm topology.
By \cite[Corollary 2.3]{Schweitzer92} the subalgebra $M_n(A^{\infty})\subseteq M_n(A)$ is closed under
the holomorphic functional calculus. It follows easily that the subalgebra $pM_n(A^{\infty})p\subseteq pM_n(A)p$ is closed
under the holomorphic functional calculus.
\end{example}

\section{Continuity of Derivations}
\label{cont of der:sec}

Throughout this section we shall assume further that $A^{\infty}$ is
closed under the smooth functional calculus. By this we mean that,
after the adjunction of a unit, for any $a\in (A^{\infty})_{\rm sa}$
and $f\in C^{\infty}(\Re)$ we have $f(a)\in A^{\infty}$.
Our goal in this section is to prove Theorem~\ref{cont of der:thm}.

\begin{example} \label{tori:example}
Let $B$ be a C*-algebra, and let $\Delta$ be a set of closed, densely defined $*$-derivations of $B$. For
$k=1, 2, {\cdots}, \infty$, define
\begin{eqnarray*}
B_k:=\{b\in B:j\le k+1,\, \delta_1, {\cdots}, \delta_j\in \Delta\Rightarrow b\in
\mbox{ the domain of }\delta_1\cdots \delta_j\}.
\end{eqnarray*}
It is routine to check that $B_k$ is a sub-$*$-algebra of $B$ and has the Fr\'echet topology determined by the seminorms
\begin{eqnarray*}
b\mapsto \pa \delta_1\cdots \delta_j b\pa, \quad j<k+1, \quad \delta_1, {\cdots}, \delta_j\in \Delta.
\end{eqnarray*}
By \cite[Lemma 3.2]{BEJ84}, if $B$ is unital, then $B_k$ is
closed under the smooth functional calculus. In particular,
$A^{\infty}_{\theta}$ is closed under the smooth functional calculus.
\end{example}

\begin{example} \label{group alg:example}
Let $G$ be a discrete group equipped with a length function $\mathnormal{l}$, that is,
a nonnegative real valued function
on $G$ such that $\mathnormal{l}(1_G)=0$,
$\mathnormal{l}(g^{-1})=\mathnormal{l}(g)$, and $\mathnormal{l}(gh)\le \mathnormal{l}(g)+
\mathnormal{l}(h)$ for all $g$ and $h$ in $G$.
Consider
the self-adjoint and closed unbounded linear operator
$D_{\mathnormal{l}}:l^2(G)\rightarrow l^2(G)$ defined by $(D_{\mathnormal{l}}\xi)(g)=\mathnormal{l}(g)\xi(g)$
for $g\in G$. Then $\delta_{\mathnormal{l}}(a)=i[D_{\mathnormal{l}}, \, a]$ defines
a closed, unbounded $*$-derivation from $\mathscr{B}(l^2(G))$ into itself.
By Example~\ref{tori:example} the intersection of the domains of $\delta^k_{\mathnormal{l}}$ for all
$k\in \Ne$ is a Fr\'echet sub-$*$-algebra of $\mathscr{B}(l^2(G))$ closed
under the smooth functional calculus. Denote by $S^{\mathnormal{l}}(G)$ the intersection
of this algebra with the reduced group algebra $C^*_r(G)$. Then $S^{\mathnormal{l}}(G)$
is a dense Fr\'echet sub-$*$-algebra of $C^*_r(G)$ closed
under the smooth functional calculus and containing $\Ce G$.
This construction is due to Connes and Moscovici \cite[page 384]{CM90}.
Recall that $G$ is said to be {\it rapidly decaying} if
there exists a length function $\mathnormal{l}$ on $G$ such that
the intersection of the domains of $D^k_{\mathnormal{l}}$
for all $k\in \Ne$, which we shall denote by $H^{\infty}_{\mathnormal{l}}(G)$,
is contained in $C^*_r(G)$ \cite{CM90, Jolissaint90}. In such
a case, it is a result of Ji \cite[Theorem 1.3]{Ji92} that
$H^{\infty}_{\mathnormal{l}}(G)$ coincides with
$S^{\mathnormal{l}}(G)$, and the Fr\'echet topology is also
induced by the seminorms $\|D^k_{\mathnormal{l}}(\cdot
)\|_{\mathnormal{l}^2}$ for $0\le k<\infty$. Consequently,
$H^{\infty}_{\mathnormal{l}}(G)$ is closed under the smooth
functional calculus if it is contained in $C^*_r(G)$.
\end{example}

\begin{question} \label{closed under sfc:question}
Is $M_m(A^{\infty})$ closed under the smooth functional calculus for all $m\in\Ne$?
\end{question}

The following theorem was proved by Sakai in the case $A^{\infty}=A$
and by Bratteli et al. in the case $A^{\infty}=A^{\infty}_{\theta}$ \cite[Corollary 5.3.C2]{BEJ84})
(cf.~also \cite[Theorem 1]{Longo79},
\cite[Theorem 3.1]{BEJ84}, and the proof of \cite[Lemma 4]{CEGJ85}).
(In fact we shall only use the result in the case $A^{\infty}=A^{\infty}_{\theta}$---in contrast to Theorem~\ref{cont of isom:thm}
which is needed in a more general setting.)

\begin{theorem} \label{cont of der:thm}
Every derivation $\delta:A^{\infty}\rightarrow A^{\infty}$ is continuous.
\end{theorem}

Theorem~\ref{cont of der:thm} is useful for determining the
derivations of various smooth (twisted) group algebras.
As an
example, let us determine the derivations of the smooth group
algebra of the $3$-dimensional discrete Heisenberg group $H_3$. This
group
is the multiplicative group
\begin{eqnarray*}
\{
\begin{pmatrix} 1 & a & c  \\ 0 &  1 & b \\ 0 & 0 & 1
\end{pmatrix}: a, b, c \in \Ze\}.
\end{eqnarray*}
It is also the universal group generated by two elements $U$ and
$V$ such that $W=VUV^{-1}U^{-1}$ is central. It is amenable
\cite[page 200]{Davidson96}.
Note that it is finitely generated.
The Fr\'echet space $H^{\infty}_l(H_3)$ does not
depend on the choice of generators if we use the word length
function corresponding to finitely many generators.
So we shall denote it by $H^{\infty}(H_3)$. Using
the word length function associated to the generators $U$ and $V$,  one
checks easily that
\begin{eqnarray} \label{smooth H3:eq}
H^{\infty}(H_3)=\{\sum_{p, q, r\in \Ze}a_{p, q, r}U^pV^qW^r\},
\end{eqnarray}
where $\{a_{p, q, r}\}$ is in the Schwartz space $\mathcal{S}(\Ze^3)$,
and the Fr\'echet topology on $H^{\infty}(H_3)$ is just the
canonical Fr\'echet topology on $\mathcal{S}(\Ze^3)$.
Note that $H_3$ has
polynomial growth, and thus
$H^{\infty}(H_3)$ is contained in $C^*(H_3)=C^*_r(H_3)$
\cite[Theorem 3.1.7]{Jolissaint90}.
By Example~\ref{group alg:example} and Theorem~\ref{cont of der:thm} we know
that every derivation of $H^{\infty}(H_3)$ into itself is continuous,
and hence is determined by the restriction on $\Ce H_3$.
Define derivations $\partial_U$ and $\partial_V$ on $\Ce H_3$ by
\begin{eqnarray*}
\partial_U(U)=U, \quad \partial_U(V)=0, \quad \partial_U(W)=0,\\
\partial_V(U)=0, \quad \partial_V(V)=V, \quad \partial_V(W)=0,
\end{eqnarray*}
and extend them continuously to $H^{\infty}(H_3)$ using (\ref{smooth H3:eq}).
We shall denote these extensions also by $\partial_U$ and $\partial_V$.
It is a result of Hadfield that every derivation $\delta:\Ce H_3\rightarrow
H^{\infty}(H_3)$ can be written uniquely
as $\delta=z_U\partial_U+z_V\partial_V+\tilde{\delta}$ for some
$z_U,\, z_V$ in the center of $H^{\infty}(H_3)$ and some $\tilde{\delta}$
as the restriction of some inner derivation of $H^{\infty}(H_3)$
\cite[Theorem 6.4]{Hadfield03}.
It is also known that the center of $H^{\infty}(H_3)$
is just the smooth algebra generated by $W$, i.e.,
$\{\sum_{r\in \Ze}a_rW^r: \{a_r\}\in \mathcal{S}(\Ze)\}$
\cite[Lemma 6.2]{Hadfield03}.
Thus we get

\begin{corollary} \label{H3:corollary}
Every derivation $\delta:H^{\infty}(H_3)\rightarrow
H^{\infty}(H_3)$ can be written uniquely
as $\delta=z_U\partial_U+z_V\partial_V+\tilde{\delta}$ for some
$z_U,\, z_V$ in the center of $H^{\infty}(H_3)$ and some inner derivation $\tilde{\delta}$ of
$H^{\infty}(H_3)$.
\end{corollary}

In view of Lemma~\ref{cont:lemma}, to prove Theorem~\ref{cont of der:thm} it suffices to
prove

\begin{lemma} \label{cont of der:lemma}
Every derivation $\delta:A^{\infty}\rightarrow A$ is continuous.
\end{lemma}

The proof of Lemma~\ref{cont of der:lemma} is similar to that of
\cite[Theorem 1]{Longo79} and \cite[Theorem 3.1]{BEJ84}. For the convenience of the reader,
we repeat the main arguments in our present setting below (in which the seminorms may not be
submultiplicative).

\begin{lemma} \label{infinite spectrum:lemma}
Let $I$ be a closed (two-sided) ideal of $A$ such that $A/I$ is
infinite-dimensional. Then there exists $b\in (A^{\infty})_{\rm sa}$
such that the image of $b$ in $A/I$ has infinite spectrum.
\end{lemma}
\begin{proof}
Without loss of generality, we may assume that $A^{\infty}$ is
unital. Consider the quotient map $\varphi:A\rightarrow A/I$. Note
that $\mathcal{A}:=\varphi(A^{\infty})$ is infinite-dimensional.
Assume that every element in $\mathcal{A}_{\rm
sa}=\varphi((A^{\infty})_{\rm sa})$ has finite spectrum (in $A/I$).
We assert that there exist nonzero projections $P_1, {\cdots}, P_m,
{\cdots}$ in $\mathcal{A}_{\rm sa}$ such that $P_jP_k=0$ for all
$j\neq k$. Assume that we have constructed $P_1, {\cdots}, P_j$ for
some $j\ge 0$ with the additional property that $Q_j\mathcal{A}Q_j$
is infinite-dimensional, where $Q_j=1-\sum^j_{s=1}P_s$. Choose an
element $b$ in $(Q_j\mathcal{A}Q_j)_{\rm sa}\setminus \Re Q_j$.
Since $b$ has finite spectrum (in $A/I$) we can find a nonzero
projection $P\in Q_j\mathcal{A}Q_j$ with $P\neq Q_j$. Since
$Q_j\mathcal{A}Q_j$ is infinite-dimensional it is easy to see that
either $P\mathcal{A}P$ or $(Q_j-P)\mathcal{A}(Q_j-P)$ has to be
infinite-dimensional. (Recall that $A/I$ is a C*-algebra.) We may
now choose $P_{j+1}$ to be one of $P$ and $Q_j-P$ in such a way that
$Q_{j+1}\mathcal{A}Q_{j+1}$ is infinite-dimensional, where
$Q_{j+1}=1-\sum^{j+1}_{s=1}P_s$.
 This finishes
the induction step.

Since $A^{\infty}$ is a Fr\'echet space we can find a complete
translation-invariant metric $d$ on $A^{\infty}$ giving the topology
of $A^{\infty}$ \cite[Corollary 13.5]{Berberian74}. Choose $b_m\in
(A^{\infty})_{\rm sa}$ such that $\varphi(b_m)=P_m$. Choose
$\lambda_m\in \Re\setminus \{0\}$ such that $d(\lambda_mb_m, 0)\le
2^{-m}$. Note that $d(\sum^m_{n=k+1}\lambda_nb_n, 0)\le 2^{-k}$ for
all $k< m$. So the series $\sum^{\infty}_{n=1}\lambda_nb_n$
converges to some $b$ in $A^{\infty}$. Since the convergence holds
in particular in $A$, it follows that
$\varphi(b)=\sum^{\infty}_{n=1}\lambda_nP_n$, the convergence being
with respect to the norm topology on $A/I$.
In particular,
we have that $\lambda_m\rightarrow 0$.
Since also $\lambda_m\neq 0$ we see that $\varphi(b)$ has infinite
spectrum (in $A/I$). This contradicts the assumption to
the contrary, which is therefore false.
In other words, $\mathcal{A}_{\rm sa}=\varphi((A^{\infty})_{\rm
sa})$ contains an element with infinite spectrum (in $A/I$).
\end{proof}

\begin{lemma} \label{finite:lemma}
Let $\delta:A^{\infty}\rightarrow A$ be a derivation. Set $\mathcal{I}=\{b\in A^{\infty}:
a\in A^{\infty}\mapsto \delta(ab)\in A \mbox{ is continuous}\}$ and denote by $I$ the closure of
$\mathcal{I}$ in $A$. Then $I$ is a closed  (two-sided) ideal of $A$, and $A/I$ is finite-dimensional.
\end{lemma}
\begin{proof} Clearly $\mathcal{I}$ is an ideal of $A^{\infty}$. So $I$ is an ideal of $A$.
Assume that $A/I$ is infinite-dimensional.
By Lemma~\ref{infinite spectrum:lemma} we can find
a self-adjoint element $b$ of $A^{\infty}$ such that the image of $b$ in $A/I$ has infinite spectrum (in $A/I$).
Choosing suitable $f_m\in C^{\infty}(\Re)$ and setting $b_m=f_m(b)$ we obtain
$\{b_m\}_{m\in\Ne}\subseteq A^{\infty}$ such that
$b^2_j\notin I$ for all $j\in \Ne$ and $b_jb_k=0$ for all $j\neq k$. We may assume that
$\pa b_j\pa , \pa \delta(b_j)\pa \le 1$ for all $j\in \Ne$.

Let $d$ be as in the proof of Lemma~\ref{infinite spectrum:lemma}.
Since $b^2_m\notin \mathcal{I}$ there exists some
 $\varepsilon_m>0$ such that for any $\varepsilon>0$ we can find
 some $a'\in A_{\infty}$ with $d(a', 0)<\varepsilon$ and $\pa
 \delta(a'b^2_m)\pa \ge \varepsilon_m$.
 The multiplication in $A^{\infty}$ is continuous by
Proposition~\ref{topological alg:prop}. Thus there exists some
$\varepsilon>0$ such that $d((m/\varepsilon_m)b'b_m, 0)\le 2^{-m}$
for any $b'\in A^{\infty}$ with $d(b', 0)<\varepsilon$.
 Take an $a'$ as above for this $\varepsilon$
 and set
 $a_m=(m/\varepsilon_m)a'$.
Then $d(a_mb_m, 0)\le 2^{-m}$ and $\pa \delta(a_mb^2_m)\pa \ge m$.
Note that $d(\sum^m_{n=k+1}a_nb_n, 0)\le 2^{-k}$ for all $k< m$. So
the series $\sum^{\infty}_{n=1}a_nb_n$ converges in $A^{\infty}$,
say to $a$. Then (with a second use of Proposition~\ref{topological
alg:prop})
$$
\pa \delta(a)\pa +\pa a\pa\ge
\pa \delta(a)b_m\pa +\pa a\delta(b_m)\pa
\ge \pa \delta(ab_m)\pa
= \pa \delta(a_mb^2_m)\pa\ge m,
$$
which is a contradiction. The assumption that $A/I$ is infinite-dimensional is therefore
not tenable. We must conclude that $A/I$ is finite-dimensional.
\end{proof}

\section{The Generic Set}
\label{gs:sec}

In this section we shall define the set $\mathcal{T}'_n$ and prove the part (1) of
Theorem~\ref{Morita equiv:thm}.

Denote by $\Der(A^{\infty}_{\theta})$ the linear space of derivations $\delta:A^{\infty}_{\theta}\rightarrow A^{\infty}_{\theta}$.
Set $\Re^n=L$. We shall think of $\Ze^n$ as the standard lattice
in $L^*$ (so that $\Hom(G, \Re)$ in \cite{BEJ84} is just our $L$),
and shall regard $\theta$ as an element of $\bigwedge^2L$. Let us write $L\otimes_{\Re}\Ce=L^{\Ce}$.
Recall that $e(t)=e^{2\pi i t}$.
One may also describe
the $C^*$-algebra $A_{\theta}$ as the universal C*-algebra generated by unitaries
$\{U_x\}_{x\in \Ze^n}$ satisfying the relations
\begin{eqnarray} \label{tori:eq}
U_xU_y=\sigma_{\theta}(x, y)U_{x+y},
\end{eqnarray}
where $\sigma_{\theta}(x, y)= e((x\cdot \theta y)/2)$. In this
description the smooth algebra $A^{\infty}_{\theta}$ becomes
$\mathcal{S}(\Ze^n, \sigma_{\theta})$, the Schwartz space
$\mathcal{S}(\Ze^n)$ equipped with the multiplication induced by
(\ref{tori:eq}). There is a canonical action of the Lie algebra
$L^{\Ce}$ as derivations of $A^{\infty}_{\theta}$, which is induced by
the canonical action of $\mathbb{T}^n$ on $A_{\theta}$ and is
given explicitly by
\begin{eqnarray*}
\delta_X(U_x)=2\pi i \left<X, x\right>U_x
\end{eqnarray*}
for $X\in L^{\Ce}$ and $x\in \Ze^n$.

\begin{notation} \label{generic:notation}
Let $e_1, {\cdots}, e_n$ be a basis of $\Ze^n$.
Denote by $\mathcal{T}'_n$ the subset of $\mathcal{T}_n$ consisting of those $\theta$'s
such that every $\delta\in \Der(A^{\infty}_{\theta})$ can be written as $\sum^n_{j=1}a_j\delta_{e_j}+\tilde{\delta}$ for some
$a_1, {\cdots}, a_n$ in the center of $A^{\infty}_{\theta}$ and some inner derivation $\tilde{\delta}$.
\end{notation}

\begin{remark} \label{generic:remark}
For any $\theta\in \mathcal{T}_n$ and $\delta\in
\Der(A^{\infty}_{\theta})$ there is at most one way of writing
$\delta$ as $\sum^n_{j=1}a_j\delta_{e_j}+\tilde{\delta}$ for some
$a_1, {\cdots}, a_n$ in the center of $A^{\infty}_{\theta}$ and
some inner derivation $\tilde{\delta}$ (see Proposition~\ref{cont of approx:prop}).
\end{remark}

One may identify $\mathcal{T}_n$ with $\Re^{\frac{n(n-1)}{2}}$ in a natural way.
We may therefore talk about Lebesgue measure on
$\mathcal{T}_n$.

\begin{proposition} \label{generic:prop}
The Lebesgue measure of $\mathcal{T}_n\setminus \mathcal{T}'_n$ is $0$.
\end{proposition}

Let $\rho_{\theta}:\Ze^n\wedge \Ze^n\rightarrow \Te$ denote the
bicharacter of $\Ze^n$ corresponding to $\theta$, i.e.,
$\rho_{\theta}(x\wedge y)=e(x\cdot \theta y)$. Recall that
$A_{\theta}$ is simple if $\rho_{\theta}$ is nondegenerate in the
sense that if $\rho_{\theta}(g\wedge h)=1$ for some $g\in \Ze^n$ and
all $h\in \Ze^n$, then $g=0$ \cite[Theorem 3.7]{Slawny72}. (In fact,
the converse is also true, though we don't need this fact here.) If
$\rho_{\theta}$ is nondegenerate,
and for every $0\neq g\in \Ze^n$ the
function $h\mapsto |\rho_{\theta}(g\wedge h)-1|^{-1} $ for
$\rho_{\theta}(g\wedge h)\neq 1$ grows at most polynomially, then
$\theta\in \mathcal{T}'_n$ \cite[page 185]{BEJ84}. So
Proposition~\ref{generic:prop} follows from

\begin{lemma} \label{measure:lemma}
Denote by $\mathcal{T}''_n$ the set of $\theta\in \mathcal{T}_n$
such that $\rho_{\theta}$ is nondegenerate
and for every $0\neq g\in \Ze^n$ the
function $h\mapsto |\rho_{\theta}(g\wedge h)-1|^{-1} $ for
$\rho_{\theta}(g\wedge h)\neq 1$ grows at most polynomially. Then
$\mathcal{T}_n\setminus \mathcal{T}''_n$ has Lebesgue measure $0$.
\end{lemma}
\begin{proof} If $\theta-\theta'\in M_n(\Ze)$, then
$\rho_{\theta}=\rho_{\theta'}$. Hence
\begin{eqnarray*}
\mathcal{T}_n\setminus \mathcal{T}''_n=\bigcup_{\eta\in M_n(\Ze)\cap
\mathcal{T}_n}(\eta\,\, +\,\, \tilde{\mathcal{T}}_n\setminus
\mathcal{T}''_n),
\end{eqnarray*}
where $\tilde{\mathcal{T}}_n$ consists of
$\theta=(\theta_{jk})\in \mathcal{T}_n$ with $0\le \theta_{jk}<1$
for all $1\le j< k\le n$. Denote Lebesgue measure on
$\mathcal{T}_n$ by $\mu$. It suffices to show that
$\mu(\tilde{\mathcal{T}}_n\setminus \mathcal{T}''_n)=0$.
If there is some polynomial $f$ in $\frac{n(n-1)}{2}$ variables such
that
\begin{eqnarray*}
1< |\sum_{1\le j<k\le n}\theta_{jk}m_{jk}-t| f(\vec{m})
\end{eqnarray*}
for all $t\in \Ze$ and $0\neq \vec{m}=(m_{jk})_{1\le j<k\le n}\in
\Ze^{\frac{n(n-1)}{2}}$, then both $\rho_{\theta}$ is nondegenerate
and
the required growth condition is satisfied---in other words, $\theta\in \mathcal{T}''_n$.
Set $F(\vec{m})=\prod_{1\le j<k\le n}(m^2_{jk}+1)^2$, and consider the set
$$
Z_{s,\vec{m}, t}=\{\theta\in \tilde{\mathcal{T}}_n:
1\ge |\sum_{1\le j<k\le n}\theta_{jk}m_{jk}-t| s F(\vec{m})\}$$
for
every $s\in \Ne$, $0\neq \vec{m}\in
\Ze^{\frac{n(n-1)}{2}}$, and $t\in \Ze$. Set $\bigcup_{\vec{m}, t}Z_{s, \vec{m},
t}=W_s$. Then every $\theta$ in $\tilde{\mathcal{T}}_n\setminus
\bigcap_{s\in \Ne}W_s$ satisfies the above condition. Therefore, it
suffices to show that $\bigcap_{s\in
\Ne}W_s$ has measure $0$. Set $\bigcup_{\substack{\vec{m}, t\\
m_{12}\neq 0}}Z_{s, \vec{m}, t}=W'_s$. Then $\mu(W_s)\le
\frac{n(n-1)}{2}\mu(W'_s)$ because of the symmetry between
the $\theta_{jk}$'s for $1\le j<k\le n$.
Integrating the characteristic function of $Z_{s,\vec{m}, t}$ over
$\theta_{12}$ first and then over the other $\theta_{jk}$'s
for $1\le j<k\le n,\, (j, k)\neq (1, 2)$, we get that
$$\mu(Z_{s,\vec{m}, t})\le
2s^{-1}(F(\vec{m}))^{-1}|m^{-1}_{12}|$$
for $m_{12}\neq 0$ and
$|t|\le |\vec{m}|:=\sum_{1\le j<k\le n}|m_{jk}|$, while
$$Z_{s,\vec{m}, t}=\ept$$
for $|t|> |\vec{m}|$. It follows that
\begin{eqnarray*}
\mu(W'_s)\le
2s^{-1}\sum_{\substack{\vec{m}\\m_{12}\neq
0}}(F(\vec{m}))^{-1}|m^{-1}_{12}|\cdot |\vec{m}|&\le &
2s^{-1}(\sum_{v\in \Ze}\frac{1}{v^2+1})^{-\frac{n(n-1)}{2}}
\\&\to &0 \quad \mbox{as } \, s\to \infty.
\end{eqnarray*}
Consequently, $\mu(\bigcap_{s\in \Ne}W_s)=0$.
\end{proof}

We shall now give two other characterizations
of $\mathcal{T}'_n$, one in Corollary~\ref{generic:corollary},
in terms of the properties of the algebra, and one in Proposition~\ref{generic 2:prop},
in terms of the number-theoretical properties of $\theta$. We need the following well-known
fact.

\begin{lemma} \label{nondegenerate:lemma}
An element $a=\sum_{h\in \Ze^n}a_hU_h$ is in the center of
$A^{\infty}_{\theta}$ if and only if the support of the
coefficients $a_h$ is contained in the subgroup $H=\{h\in \Ze^n:
\rho_{\theta}(g\wedge h)=1 \mbox{ for all } g\in \Ze^n\}$. In
particular, the center of $A^{\infty}_{\theta}$ is $\Ce$ if and
only if $\rho_{\theta}$ is nondegenerate in the sense that
$H=\{0\}$.
\end{lemma}
\begin{proof} The element $a$ is in the center
of $A^{\infty}_{\theta}$ exactly if $U_ga=aU_g$ for all
$g\in \Ze^n$. Note that $U_ga(U_g)^{-1}=\sum_{h\in \Ze^n}a_h\rho_{\theta}(g\wedge h)U_h$.
Consequently, $U_ga=aU_g$ for all $g\in \Ze^n$ if and only if $a_h=0$ for
all $h\in \Ze^n\setminus H$.
\end{proof}

Recall that the topology considered on $A^{\infty}_{\theta}$ will
always be the Fr\'echet topology, unless otherwise specified.

\begin{definition} \label{ader:def}
Let us say that a derivation $\delta\in \Der(A^{\infty}_{\theta})$
is {\it approximately inner} if there is a sequence
$\{a_m\}_{m\in\Ne}\subseteq A^{\infty}_{\theta}$ such that $[a_m,
a]\to\delta(a)$ (in the Fr\'echet topology) for every $a\in
A^{\infty}_{\theta}$. Denote by $\ADer(A^{\infty}_{\theta})$ the
linear space of approximately inner derivations of
$A^{\infty}_{\theta}$.
\end{definition}

By \cite[Corollary 5.3.D2]{BEJ84}, every $\delta \in
\Der(A^{\infty}_{\theta})$ can be written uniquely as
$\sum^n_{j=1}a_j\delta_{e_j}+\tilde{\delta}$ for some $a_1,
{\cdots}, a_n$ in the center of $A^{\infty}_{\theta}$ and some
$\tilde{\delta}\in \Der(A^{\infty}_{\theta})$ such that there is a
sequence $\{b_m\}_{m\in\Ne}\subseteq A^{\infty}_{\theta}$ with $\pa
[b_m, a]-\tilde{\delta}(a)\pa \to 0$ for every $a\in
A^{\infty}_{\theta}$. In fact, we can require $[b_m, a]\to
\tilde{\delta}(a)$ in the Fr\'echet topology:

\begin{proposition} \label{cont of approx:prop}
Every $\delta \in \Der(A^{\infty}_{\theta})$ can be written
uniquely as $\sum^n_{j=1}a_j\delta_{e_j}+\tilde{\delta}$ for some $a_1, {\cdots}, a_n$ in the center of
$A^{\infty}_{\theta}$ and some $\tilde{\delta}\in \ADer(A^{\infty}_{\theta})$.
The bicharacter $\rho_{\theta}$ is nondegenerate
if and only if every
$\delta\in \Der(A^{\infty}_{\theta})$ can
be written uniquely as $\delta_X+\tilde{\delta}$ for some $X\in L^{\Ce}$ and
some $\tilde{\delta}\in \ADer(A^{\infty}_{\theta})$.
\end{proposition}

Denote by $A^F_{\theta}$ the linear span of $\{U_x\}_{x\in \Ze^n}$.
This is a dense sub-$*$-algebra of $A_{\theta}$. In the proof of
\cite[Corollary 5.3.D2]{BEJ84}, which itself is based on the proof
of \cite[Theorem 2.1]{BEJ84}, one sees easily that in the present
case actually the sequence $\{b_m\}_{m\in\Ne}$ can be chosen in such
a way that $[b_m, a]\to \tilde{\delta}(a)$ in the Fr\'echet topology
for every $a\in A^F_{\theta}$. Thus
Proposition~\ref{cont of approx:prop} follows from
the following lemma and Lemma~\ref{nondegenerate:lemma}.

\begin{lemma} \label{conv of der:lemma}
Let $\delta, \delta_1, \delta_2, {\cdots}\in \Der(A^{\infty}_{\theta})$
be such that $\delta_m(a)\to \delta(a)$ for every $a\in A^F_{\theta}\subseteq A^{\infty}_{\theta}$.
Then  $\delta_m(a)\to \delta(a)$ for every $a\in A^{\infty}_{\theta}$.
\end{lemma}
\begin{proof}
The proof is similar to that of \cite[Corollary 3.3.4]{BEJ84}.
Let $\vec{j}=(j_1, {\cdots}, j_k)$ with $1\le j_1, {\cdots}, j_k\le n$ and $k\ge 0$.
We say that $\vec{w}=(w_1, {\cdots}, w_s)$ is a subtuple of $\vec{j}$ if $s\le k$ and there is
a strictly increasing map $f:\{1, {\cdots}, s\}\rightarrow \{1, {\cdots}, k\}$ such that
$w_m=j_{f(m)}$.
Set $\pa a\pa_{\vec{j}}=\sup \pa \delta_{w_s}\cdots \delta_{w_1}(a)\pa$, where the $\sup$ runs
over all subtuples $\vec{w}$ of $\vec{j}$,
for $a\in A^{\infty}_{\theta}$.
It suffices to show that $\pa \delta(a)-\delta_m(a)\pa_{\vec{j}}\to 0$ for every $a\in A^{\infty}_{\theta}$
and $\vec{j}$. For each $g=(q_1, {\cdots}, q_n)\in \Ze^n$ set $|g|=\sum^s_{s=1}|q_s|$.
Set $M=\sup\{\pa \delta_m(U^s_t)\pa_{\vec{j}}:t=1, {\cdots}, n; s=\pm 1; m=1, 2, \cdots\}$.
Using the derivation property of $\delta_m$ we have
$\pa \delta_m(U^{q_1}_1\cdots U^{q_n}_n)\pa_{\vec{j}}\le (2\pi)^kM|g|^{k+1}$ for every
$g=(q_1, {\cdots}, q_n)\in \Ze^n$.
For any finite subset $Z\subseteq \Ze^n$ and $a=\sum_{g\in \Ze^n}a_gU_g\in A^{\infty}_{\theta}$
set $a_Z=\sum_{g\in Z}a_gU_g\in A^F_{\theta}$. By Theorem~\ref{cont of der:thm} every derivation on
$A^{\infty}_{\theta}$ is continuous. Thus $\pa \delta(a-a_Z)\pa_{\vec{j}}\to 0$ as $Z$
goes to $\Ze^n$. Also
$\pa \delta_m(a-a_Z)\pa_{\vec{j}}\le
\sum_{g\in \Ze^n\setminus Z}(2\pi)^kM|g|^{k+1}|a_g|$. So $\pa \delta_m(a-a_Z)\pa_{\vec{j}}\to 0$ uniformly
as $Z$ goes to $\Ze^n$. By assumption $\pa \delta(a_Z)-\delta_m(a_Z)\pa_{\vec{j}}\to 0$ as $m\to \infty$.
Thus $\pa \delta(a)-\delta_m(a)\pa_{\vec{j}}\to 0$ as $m\to \infty$.
\end{proof}

\begin{corollary} \label{approx inner:coro}
A derivation $\delta\in \Der(A^{\infty}_{\theta})$ is
approximately inner if and only if there is a sequence
$\{b_m\}_{m\in\Ne}\subseteq A^{\infty}_{\theta}$ with $\pa [b_m,
a]-\delta(a)\pa \to 0$ for every $a\in A^{\infty}_{\theta}$.
\end{corollary}

Combining Lemma~\ref{nondegenerate:lemma} and Proposition~\ref{cont of approx:prop}
we get

\begin{corollary} \label{generic:corollary}
The set $\mathcal{T}'_n$ consists of those $\theta$'s such that
every $\delta\in \ADer(A^{\infty}_{\theta})$ is
inner.
\end{corollary}

Set $F(g)=\max_{1\le j\le n}|\rho_{\theta}(g\wedge e_j)-1|$ for
$g\in \Ze^n$. Then $F$ vanishes exactly on the subgroup $H$ in Lemma~\ref{nondegenerate:lemma}.
Denote by $F^{-1}$ the function on $\Ze^n$ taking on the value $F(g)^{-1}$ at $g\notin H$ and the value
$0$ at $g\in H$.

\begin{proposition} \label{generic 2:prop}
The set $\mathcal{T}'_n$ consists of those $\theta$'s such that
the function
$F^{-1}$ grows at most polynomially.
\end{proposition}
\begin{proof}
In view of
Corollary~\ref{generic:corollary}
it suffices to show that
every $\delta\in \ADer(A^{\infty}_{\theta})$ is
inner if and only if the function
$F^{-1}$ grows at most polynomially.

By Proposition~\ref{cont of approx:prop},
Theorem~\ref{cont of der:thm}, \cite[Corollary 5.3.E2]{BEJ84},
and the proof of \cite[Theorem 5.1]{BEJ84}
(see also the first paragraph of the proof of  \cite[Theorem 2.1]{BEJ84})
the derivations $\delta\in \ADer(A^{\infty}_{\theta})$ are in bijective correspondence with
those $\Ce$-valued functions $Q$ on $\Ze^n$ such that
$Q$ vanishes on $H$ and the function $c_h:g\mapsto Q(g)(\rho_{\theta}(g\wedge h)-1)$ on $\Ze^n$ is in
the Schwarz space $\mathcal{S}(\Ze^n)$
for every $h\in \Ze^n$. Actually $\delta(U_h)=\sum_{g}c_h(g)U_hU_g$.
Furthermore, $\delta$ is inner
if and only if $Q\in \mathcal{S}(\Ze^n)$.
In this case, $\delta(\cdot)=[\sum_gQ(g)U_g, \cdot]$.
Clearly, $c_h\in \mathcal{S}(\Ze^n)$ for every $h\in \Ze^n$ if and only if
$c_{e_j}\in \mathcal{S}(\Ze^n)$ for every $1\le j\le n$, and if and only if
the function $g\mapsto Q(g)F(g)$ is in $\mathcal{S}(\Ze^n)$.
In other words, the derivations $\delta\in \ADer(A^{\infty}_{\theta})$ are in bijective correspondence with
those $Q:\Ze^n\rightarrow \Ce$ such that $Q$ vanishes on $H$ and
the function $g\mapsto Q(g)F(g)$ is in $\mathcal{S}(\Ze^n)$.
Therefore, every $\delta\in \ADer(A^{\infty}_{\theta})$ is
inner if and only if the pointwise multiplication by $F^{-1}$ sends $\mathcal{S}(\Ze^n)$ into itself.
Using the closed graph theorem \cite[Corollary 48.6]{Berberian74} it is easy to see that if
the pointwise multiplication by $F^{-1}$ sends $\mathcal{S}(\Ze^n)$ into itself, then
this map is continuous and hence $F^{-1}$ grows at most polynomially. Conversely,
if $F^{-1}$ grows at most polynomially, then obviously
the pointwise multiplication by $F^{-1}$ sends $\mathcal{S}(\Ze^n)$ into itself.
Therefore every $\delta\in \ADer(A^{\infty}_{\theta})$ is
inner if and only if $F^{-1}$
grows at most polynomially.
\end{proof}

Denote by $\mathcal{T}^{\flat}_n$ the subset of $\cT_n$ consisting
of $\theta$'s such that $\rho_{\theta}$ is nondegenerate. Let us
indicate how to deduce the weaker form of the part (1) of
Theorem~\ref{Morita equiv:thm} with $\mathcal{T}'_n$ replaced by
$\mathcal{T}'_n\cap \mathcal{T}^{\flat}_n$ from \cite{Nest88} using
the topological or algebraic Hochschild cohomology of
$A^{\infty}_{\theta}$. Recall that if two unital algebras are Morita
equivalent, then their algebraic Hochschild cohomology is isomorphic
\cite{Loday}. By Theorem~\ref{cont of isom:thm} and Example~\ref{cut
down:example} if two unital smooth algebras are Morita equivalent,
then their topological Hochschild cohomology is also isomorphic.
Nest calculated the topological Hochschild cohomology ${\rm
H}^*_{\rm top}(A^{\infty}_{\theta}, (A^{\infty}_{\theta})^*_{\rm
top})$ of (the Fr\'echet algebra) $A^{\infty}_{\theta}$ with
coefficients in the topological dual $(A^{\infty}_{\theta})^*_{\rm
top}$
in \cite[Theorem 4.1]{Nest88} (the $2$-dimensional case was
calculated earlier by Connes in \cite{Connes85}).
From \cite[Theorem 4.1]{Nest88} it is
easy to see that $\rho_{\theta}$ is nondegenerate and $\theta$
satisfies the condition in Proposition~\ref{generic 2:prop}
if and only if  ${\rm H}^*_{\rm top}(A^{\infty}_{\theta},
(A^{\infty}_{\theta})^*_{\rm top})$ is finite-dimensional in every
degree. Using the simple projective resolution of
$A^{\infty}_{\theta}$ as an $A^{\infty}_{\theta}$-bimodule in
\cite[Section 3]{Nest88} one also finds that this happens if and
only if the algebraic Hochschild cohomology ${\rm H}^*_{\rm
alg}(A^{\infty}_{\theta}, (A^{\infty}_{\theta})^*_{\rm alg})$ is
finite-dimensional in every degree. Thus the above weak form of the
part (1) of Theorem~\ref{Morita equiv:thm} follows from considering
either the topological or algebraic Hochschild cohomology of
$A^{\infty}_{\theta}$.

In the $2$-dimensional case, when $\rho_{\theta}$ is nondegenerate,
the condition in Proposition~\ref{generic 2:prop} was called a {\it
diophantine condition} by Connes \cite[page 349]{Connes85}.

In \cite{Boca97} Boca introduced a certain subset of $\mathcal{T}_n$
the complement of which has Lebesgue measure $0$ and which is also
described number theoretically. His set is contained in
$\mathcal{T}^{\flat}_n$. We do not know whether his set is the same
as $\mathcal{T}'_n\cap \mathcal{T}^{\flat}_n$ or not.

To prove the part (1) of Theorem~\ref{Morita equiv:thm} we start
with some general facts about the comparison of derivation spaces
for Morita equivalent algebras. Let $\cA$ be a unital algebra.
Let $E$ be a finitely generated projective right $\cA$-module and
set $\End(E_{\cA})=\cB$.
If we take an isomorphism of right $\cA$-modules $E\rightarrow
p(k\cA)$ for some projection $p\in M_k(\cA)$, where $k\cA$ is the
direct sum of $k$ copies of $\cA$ as right $\cA$-modules with
vectors written as columns, then we have an induced isomorphism
$\cB\rightarrow pM_k(\cA)p$.

Let $\delta\in \Der(\cA)$. Recall that a {\it connection}
\cite{Connes80} for $(E_{\cA}, \delta)$ is a linear map
$\nabla:E\rightarrow E$ satisfying the Leibnitz rule
\begin{eqnarray}
\nabla(fa)=\nabla(f)a+f\delta(a)
\end{eqnarray}
for all $f\in E$ and $a\in \cA$. Let us say that a pair $(\delta',
\delta)\in \Der(\cB)\times \Der(\cA)$ is {\it compatible} if there
is a linear map $\nabla:E\rightarrow E$ which is a connection for
both $({}_{\cB}E, \,\delta')$ and $(E_{\cA},\, \delta)$. One checks
easily that for every $\delta\in \Der(\cA)$ there exists $\delta'\in
\Der(\cB)$ such that the pair $(\delta',\, \delta)$ is compatible,
and $\delta'$ is unique up to adding an inner derivation.
Explicitly, identifying $E$ and $\cB$ with $p(k\cA)$ and
$pM_k(\cA)p$ respectively as above, and extending $\delta$ to $k\cA$
and $M_k(\cA)$ componentwise, one may choose $\nabla$ and $\delta'$
as defined by $\nabla(u)=p(\delta(u))$ for $u\in p(k\cA)$ and
$\delta'(b)=p\delta(b)p$ for $b\in pM_k(\cA)p$ respectively.

\begin{lemma} \label{inner:lemma}
If $\delta$ is inner and the pair $(\delta',\, \delta)$ is compatible, then $\delta'$ is also inner.
\end{lemma}
\begin{proof}
Say $\delta(\cdot)=[\cdot, \, a]$ for some $a\in \cA$. The pair
$(0,\, \delta)$ is compatible with respect to the connection
$\nabla_a: f\mapsto fa$. Say the pair $(\delta',\, \delta)$ is
compatible with respect to the connection $\nabla$. Then the pair
$(\delta',\, 0)$ is compatible with respect to the connection
$\nabla-\nabla_a$. Therefore, $\delta'$ is inner.
\end{proof}

Assume further that $\cA=A^{\infty}$ is equipped with a smooth
topology. By Example~\ref{cut down:example} and Theorem~\ref{cont of
isom:thm} we know that $\cB=B^{\infty}$ also admits a unique smooth
topology. Thus, the above isomorphism $B^{\infty}\rightarrow
pM_k(A^{\infty})p$ is a homeomorphism.

\begin{lemma} \label{app inner:lemma}
If $\delta\in \ADer(A^{\infty})$ and the pair $(\delta', \delta)$ is compatible,
then $\delta'\in \ADer(B^{\infty})$.
\end{lemma}
\begin{proof} We may assume that $E=p(kA^{\infty})$ and $B^{\infty}=pM_k(A^{\infty})p$
for some projection $p\in M_k(A^{\infty})$. Choose a sequence of inner derivations
$\delta_m\in \Der(A^{\infty})$ such that
$\delta_m(a)\to \delta(a)$ for every $a\in A^{\infty}$. Extend $\delta_m$ and $\delta$ to
$M_k(A^{\infty})$ componentwise.
Consider the maps $\tilde{\delta}:b\to p\delta(b)p$ and $\tilde{\delta}_m:
b\to p\delta_m(b)p$ on $B^{\infty}$. Then $\tilde{\delta}$ and $\tilde{\delta}_m$
are derivations of $B^{\infty}$. The pair $(\tilde{\delta}, \delta)$ is compatible, with
respect to the Grassmann connection $\nabla: u\mapsto p(\delta(u))$. Similarly,
the pair $(\tilde{\delta}_m, \delta_m)$ is compatible,  with
respect to the connection $\nabla_m: u\mapsto p(\delta_m(u))$.
Notice that
$\tilde{\delta}_m(b)\to \tilde{\delta}(b)$ for every $b\in B^{\infty}$.
Lemma~\ref{app inner:lemma} now follows from Lemma~\ref{inner:lemma}.
\end{proof}

There are several equivalent ways of defining Morita equivalence of
algebras (see for instance \cite[Section 22]{AF74}). Recall that a
right $\mathcal{A}$-module $E_{\mathcal{A}}$ of a unital algebra
$\mathcal{A}$ is a {\it generator} if $\mathcal{A}_{\mathcal{A}}$ is
a direct summand of $rE_{\mathcal{A}}$ for some $r\in \Ne$. We shall
say that two unital algebras $\mathcal{B}$ and $\mathcal{A}$ are
{\it Morita equivalent} if there exists a bimodule
${}_{\mathcal{B}}E_{\mathcal{A}}$---a Morita equivalence
bimodule---such that $E_{\mathcal{A}}$ and ${}_{\mathcal{B}}E$ are
finitely generated projective modules and also generators, and,
furthermore, $\mathcal{B}=\End(E_{\mathcal{A}}),\,
\mathcal{A}=\End({}_{\mathcal{B}}E)$ \cite[Theorem 22.2]{AF74}.

Now the part (1) of Theorem~\ref{Morita
equiv:thm} follows from Corollary~\ref{generic:corollary},
Lemma~\ref{inner:lemma} and Lemma~\ref{app inner:lemma}.

The above proof employs the Fr\'echet topology on
$A^{\infty}_{\theta}$. We give below a more algebraic proof. We are
grateful to Ryszard Nest for suggesting using the Morita invariance
of the module structure of $\rH^1(\cA, \cA)$ over the center of
$\cA$
 for a unital algebra $\cA$.

The part (1) of Theorem~\ref{Morita equiv:thm} follows directly from
two facts. Denote by $Z(\cA)$ the center of a unital algebra $\cA$.
Given a Morita equivalence bimodule $_{\cB}E_{\cA}$ between two
unital algebras $\cB$ and $\cA$, we may identify both $Z(\cB)$ and
$Z(\cA)$ with $\End(E_{\cA})\cap \End(_{\cB}E)$ inside
$\Hom_{\Ce}(E)$. Note that $\Der(\cA)$ has a natural $Z(\cA)$-module
structure given by $(a\delta)(x)=a(\delta(x))$ for all $a\in
Z(\cA)$, $x\in \cA$, and $\delta\in \Der(\cA)$. Clearly the space of
inner derivations is a submodule. Denote by $\Out(\cA)$ the quotient
$Z(\cA)$-module. By Lemma~\ref{inner:lemma} we have a natural linear
isomorphism $\Out(\cA)\rightarrow \Out(\cB)$. The first fact we need
is that this linear isomorphism is (clearly) an isomorphism of
$Z(\cA)$($=Z(\cB)$)-modules. The second fact is that $\theta\in
\mathcal{T}_n$ if and only if $\Out(A^{\infty}_{\theta})$ is
generated by $n$ elements as a $Z(A^{\infty}_{\theta})$-module. This
follows from the decomposition of $\Der(A^{\infty}_{\theta})$ quoted
after Definition~\ref{ader:def}.

One may also deduce the second fact as follows. Recall that the
first (algebraic) Hochschild cohomology $\rH^1(\cA, \cA)$ of $\cA$
with coefficients in $\cA$ is exactly $\Out(\cA)$ \cite[page
38]{Loday}. Using the simple projective resolution of
$A^{\infty}_{\theta}$ as an $A^{\infty}_{\theta}$-bimodule in
\cite[Section 3]{Nest88} one can calculate
$\rH^*(A^{\infty}_{\theta}, A^{\infty}_{\theta})$ and find that
$\theta$ satisfies the condition in Proposition~\ref{generic 2:prop}
if and only if $\rH^1(A^{\infty}_{\theta}, A^{\infty}_{\theta})$ is
generated by $n$ elements as a $Z(A^{\infty}_{\theta})$-module.

Combining Lemma~\ref{nondegenerate:lemma}, Proposition~\ref{cont of approx:prop}, and
 Lemma~\ref{app inner:lemma}
we also get

\begin{proposition} \label{linear map:prop}
Suppose that $\rho_{\theta}$ and $\rho_{\theta'}$ are nondegenerate,
and that $E$ is a finitely generated
projective right $A^{\infty}_{\theta}$-module with
$\End(E_{A^{\infty}_{\theta}})=A^{\infty}_{\theta'}$.
Then there is a unique linear map $\varphi: L^{\Ce}\rightarrow L^{\Ce}$ such that
for any $X\in L^{\Ce}$ and $\tilde{\delta}\in \ADer(A^{\infty}_{\theta})$ there
exists some $\tilde{\delta}'\in \ADer(A^{\infty}_{\theta'})$ such that
the pair $(\delta_{\varphi(X)}+\tilde{\delta}', \delta_X+\tilde{\delta})$ is compatible.
If, furthermore, the bimodule ${}_{A^{\infty}_{\theta'}}E_{A^{\infty}_{\theta}}$ is a Morita equivalence bimodule,
then $\varphi$ is an isomorphism.
\end{proposition}

\section{Proof of Theorem~\ref{Morita equiv:thm}}
\label{proof:sec}

In this section we prove the part (2) of Theorem~\ref{Morita equiv:thm}.

We recall first the theory of curvature introduced by Connes in
\cite{Connes80}. Let $E$ be a finitely generated projective right
$A^{\infty}_{\theta}$-module. If $X\in L^{\Ce}\mapsto \nabla_X\in
\Hom_{\Ce}(E)$ is a linear map such that $\nabla_X$ is a connection
of $(E_{A^{\infty}_{\theta}}, \delta_X)$ for every $X\in L^{\Ce}$,
one may consider the {\it curvature} $[\nabla_X, \nabla_Y]$ which is
easily seen to be in $\End(E_{A^{\infty}_{\theta}})$. We say that
$\nabla$ has {\it constant curvature} if $[\nabla_X, \nabla_Y]\in
\Ce$($=\Ce \cdot \id_{E}$) for all $X, Y\in L^{\Ce}$.

Since complete Morita equivalence in the sense of Schwarz
\cite{Schwarz98} explicitly implies Morita equivalence (see \cite[Subsection 2.1]{EL2}),
the ``if'' part of
the statement follows from \cite[Theorem 1.2]{LiMorita} (which deals
with complete Morita equivalence).

Recall that $\mathcal{T}^{\flat}_n$ is the subset of $\cT_n$
consisting of $\theta$'s such that $\rho_{\theta}$ is nondegenerate.
To prove the ``only if'' part of the statement for all $n$, we shall
reduce it first to the case of matrices in
$\mathcal{T}'_n\cap \mathcal{T}^{\flat}_n$. For this purpose, we need
the following lemma, in which we consider also $\cT_0$ for convenience.

\begin{lemma} \label{simple quotient:lemma}
Suppose that $\theta\in \cT_n$ is of the form
\begin{eqnarray*}
\theta=\begin{pmatrix} 0 & 0 \\ 0 &  \tilde{\theta} \end{pmatrix},
\end{eqnarray*}
where $\tilde{\theta}$ belongs to $\cT_k$ for some $0\le k\le n$
and $\rho_{\tilde{\theta}}$ is nondegenerate. Then for any maximal two-sided ideal $I$
of $A^{\infty}_{\theta}$, $A^{\infty}_{\theta}/I$ is
isomorphic to $A^{\infty}_{\tilde{\theta}}$.
\end{lemma}
\begin{proof}
Since $A^{\infty}_{\theta}$ is closed under the smooth functional
calculus (Example~\ref{tori:example}), by Theorem 13 of \cite{KS}
and the remark following it there is a bijective correspondence
between the lattice of two-sided ideals of $A^{\infty}_{\theta}$
closed with respect to the relative C*-algebra topology and the
lattice of closed two-sided ideals of $A_{\theta}$ which in one
direction consists in taking the intersection of an ideal with
$A^{\infty}_{\theta}$ and in the other direction in taking the
closure in $A_{\theta}$. Since $I$ is maximal, it is closed in the
relative C*-algebra topology. It follows that the closed two-sided
ideal $K$ of $A_{\theta}$ corresponding to $I$ is maximal. Note that
$A_{\theta}=C(\Te^{n-k})\otimes A_{\tilde{\theta}}$ and that the
center of $A_{\theta}$ is $C(\Te^{n-k})$. Since
$\rho_{\tilde{\theta}}$ is nondegenerate, $A_{\tilde{\theta}}$ is
simple \cite[Theorem 3.7]{Slawny72}. It follows that $K$ is equal to
the kernel of the homomorphism $A_{\theta}\rightarrow
A_{\tilde{\theta}}$ given by the evaluation at some point of
$\Te^{n-k}$. Then we may identify $A_{\theta}/K$ with
$A_{\tilde{\theta}}$. It follows that $A^{\infty}_{\theta}/I$, which
is contained in $A_{\theta}/K$, is just
$A^{\infty}_{\tilde{\theta}}$.
\end{proof}

Suppose that $\theta',\, \theta\in \mathcal{T}'_n$, and  that
$A^{\infty}_{\theta'}$ and $A^{\infty}_{\theta}$ are Morita
equivalent, By \cite[Proposition 3.3]{EL2}, we can find $\theta'_1,
\, \theta_1 \in \cT_n$ such that $\theta'$ and $\theta$ are in the
same orbit of the $\SO(n, n|\Ze)$ action as $\theta'_1$ and
$\theta_1$ respectively and such that
\begin{eqnarray*}
\theta'_1=\begin{pmatrix} 0 & 0 \\ 0 &  \tilde{\theta'} \end{pmatrix},\quad \theta_1=\begin{pmatrix} 0 & 0 \\ 0 &  \tilde{\theta} \end{pmatrix},
\end{eqnarray*}
where $\tilde{\theta'}$ and $\tilde{\theta}$ belong to $\cT_{k'}$
and $\cT_{k}$ respectively for some $0\le k',\, k\le n$ and
$\rho_{\tilde{\theta'}}$ and $\rho_{\tilde{\theta}}$ are
nondegenerate. By \cite[Theorem 1.2]{LiMorita},
$A^{\infty}_{\theta'}$ and $A^{\infty}_{\theta}$ are completely
Morita equivalent to $A^{\infty}_{\theta'_1}$ and
$A^{\infty}_{\theta_1}$ respectively. Then
$A^{\infty}_{\theta'_1}$ and $A^{\infty}_{\theta_1}$ are Morita
equivalent. Since Morita equivalence between unital algebras (or
rings) preserves the center \cite[Proposition 21.10]{AF74}, by
Lemma~\ref{nondegenerate:lemma} we have $n-k'=n-k$. Therefore,
$k'=k$. There is a natural bijection between the lattices of
two-sided ideals of Morita equivalent unital algebras (or rings),
and the corresponding quotient algebras are also Morita equivalent
\cite[Proposition 21.11]{AF74}. It follows from Lemma~\ref{simple
quotient:lemma} that $A^{\infty}_{\tilde{\theta'}}$ and
$A^{\infty}_{\tilde{\theta}}$ are Morita equivalent. By the part
(1) of Theorem~\ref{Morita equiv:thm}, $\theta'_1$ and $\theta_1$
are also in $\cT'_n$. From Proposition~\ref{generic 2:prop} we see
that $\tilde{\theta'}$ and $\tilde{\theta}$ are both in $\cT'_k$.
If the ``only if'' part of the part (2) of Theorem~\ref{Morita
equiv:thm} holds for all $n$ with $\cT'_n$ replaced by
$\cT'_n\cap\cT^{\flat}_n$, then we can conclude that
$\tilde{\theta'}$ and $\tilde{\theta}$ are in the same orbit of
the $\SO(k, k|\Ze)$ action. It follows that $\theta'_1$ and
$\theta_1$ are in the same orbit of the $\SO(n, n|\Ze)$ action.
Consequently, $\theta'$ and $\theta$ are in the same orbit of the
$\SO(n, n|\Ze)$ action.

Now what remains is to prove
the ``only if'' part of the
statement with $\cT'_n$ replaced by $\cT'_n\cap\cT^{\flat}_n$. This
follows from Theorems~\ref{complete me to orbit:thm} and
\ref{me to complete:thm} below.

\begin{theorem} \label{complete me to orbit:thm}
Suppose that $A^{\infty}_{\theta'}$ and $A^{\infty}_{\theta}$ are Morita equivalent with respect
to a Morita equivalence bimodule
$E={}_{A^{\infty}_{\theta'}}E_{A^{\infty}_{\theta}}$ with the following property: there are a $\Ce$-linear isomorphism
$\varphi: L^{\Ce}\rightarrow L^{\Ce}$
and a linear map $\nabla$ from $L^{\Ce}$
to $\Hom_{\Ce}(E)$,
with $\nabla_X$  a connection for
both $(E_{A^{\infty}_{\theta}}, \delta_X)$ and $({}_{A^{\infty}_{\theta'}}E, \delta_{\varphi(X)})$
for every $X\in L^{\Ce}$, such that
$\nabla$ has constant curvature (in both $\End(E_{A^{\infty}_{\theta}})$ and $\End({}_{A^{\infty}_{\theta'}}E)$).
Then $\theta'$ and $\theta$ are in the same orbit of the $\SO(n, n|\, \Ze)$ action.
\end{theorem}

\begin{theorem} \label{me to complete:thm}
Suppose that $\theta',\, \theta\in \mathcal{T}'_n\cap \cT^{\flat}_n$,
and  that
$A^{\infty}_{\theta'}$ and $A^{\infty}_{\theta}$ are Morita equivalent,
and let ${}_{A^{\infty}_{\theta'}}E_{A^{\infty}_{\theta}}$
be a Morita equivalence bimodule. Then there exist $\varphi$ and
$\nabla$ satisfying the conditions of Theorem~\ref{complete me to orbit:thm}. (N.B: $E$ and $\nabla$
are not necessarily Hermitian---see the proof of Theorem~\ref{complete me to orbit:thm} below.)
\end{theorem}
\begin{proof}[Proof of Theorem~\ref{me to complete:thm}]
Let $\varphi$ be as in Proposition~\ref{linear map:prop}. Let $e_1,
{\cdots}, e_n$ be a basis of $L^{\Ce}$. Then for each $1\le k\le n$
there is some $\delta'_k\in \ADer(A^{\infty}_{\theta'})$ such that
$(\delta_{\varphi(e_k)}+\delta'_k, \delta_{e_k})$ is compatible.
Since $\theta'\in \mathcal{T}'_n$, by
Corollary~\ref{generic:corollary} the derivation $\delta'_k$ is
inner. Then $(\delta_{\varphi(e_k)}, \delta_{e_k})$ is compatible.
Let $\nabla_{e_k}$ be a connection for both
$({}_{A^{\infty}_{\theta'}}E, \delta_{\varphi(e_k)})$ and
$(E_{A^{\infty}_{\theta}}, \delta_{e_k})$. Let $ \nabla_{\sum_{1\le
k\le n}\lambda_ke_k}=\sum_{1\le k\le n}\lambda_k\nabla_{e_k}$ for
all $\lambda_1, {\cdots}, \lambda_n\in \Ce$. Then for every $X\in
L^{\Ce}$, $\nabla_X$ is a connection for both
$({}_{A^{\infty}_{\theta'}}E, \delta_{\varphi(X)})$ and
$(E_{A^{\infty}_{\theta}}, \delta_X)$. Consequently, $[\nabla_X,\,
\nabla_Y]\in \Hom_{\Ce}(E)$ is in both $A^{\infty}_{\theta'}$ and
$A^{\infty}_{\theta}$. Therefore, $[\nabla_X,\, \nabla_Y]$ lies in
the center of $A^{\infty}_{\theta}$. By
Lemma~\ref{nondegenerate:lemma} the center of $A^{\infty}_{\theta}$
is $\Ce$. Consequently, $[\nabla_X, \nabla_Y] \in \Ce$.
\end{proof}

Theorem~\ref{complete me to orbit:thm} is an extension of Schwarz's result of \cite[Section 5]{Schwarz98}, in which he proved
Theorem~\ref{complete me to orbit:thm} under the addition hypotheses that $E$ is a Hilbert bimodule,
$\varphi$ maps $L$ to $L$,
and $\nabla$ is a Hermitian connection. We shall essentially follow Schwarz's argument.
In order to show that his argument still works without these additional hypotheses, we have to make some
preparations.

Let $E_{A^{\infty}_{\theta}}$ be a finitely generated projective right $A^{\infty}_{\theta}$-module.
Let $\tau$ be a trace on $A^{\infty}_{\theta}$.
Denote by $\tau$ also the unique extension of $\tau$ to a trace of $M_k(A^{\infty}_{\theta})$ for
each $k\in \Ne$.
Choosing an isomorphism $E\rightarrow p(kA^{\infty}_{\theta})$
for some idempotent $p\in M_k(A^{\infty}_{\theta})$,
we get a trace
$\tau'$ on $\End(E_{A^{\infty}_{\theta}})$, via the natural isomorphism of this algebra with
$pM_k(A^{\infty}_{\theta})p$. It is not difficult to see that $\tau'$ does not depend on the
choice of $p$ or the isomorphism $E\rightarrow p(kA^{\infty}_{\theta})$.

\begin{lemma} \label{trace and der:lemma}
Let $\delta\in \Der(A^{\infty}_{\theta})$ be such that $\tau \circ \delta=0$. Then
$\tau'\circ \delta'=0$ for any $\delta'\in \Der(\End(E_{A^{\infty}_{\theta}}))$
such that the pair $(\delta', \delta)$ is compatible.
\end{lemma}
\begin{proof} We may assume that $E=p(kA^{\infty}_{\theta})$ and $\End(E_{A^{\infty}_{\theta}})=pM_k(A^{\infty}_{\theta})p$
for some idempotent $p\in M_k(A^{\infty}_{\theta})$.
Extend $\delta$ to
$M_k(A^{\infty}_{\theta})$ componentwise.
Let $\tilde{\delta}$ be the derivation of $pM_k(A^{\infty}_{\theta})p$ defined
as $\tilde{\delta}(b)=p\delta(b)p$. Then the pair $(\tilde{\delta}, \delta)$ is compatible.
For each $b\in pM_k(A^{\infty}_{\theta})p$ pick $u_j, v_j\in kA^{\infty}_{\theta}$ for $1\le j\le k$ such that
$b=\sum_{1\le j\le k}u_jv^t_j$. Then $b=pbp=\sum_{1\le j\le k}(pu_j)(v^t_jp)$, so we may assume that
$pu_j=u_j$ and $v^t_jp=v^t_j$. Now
\begin{eqnarray*}
\tau(\tilde{\delta}(b))&=&\tau(\sum_j(p\delta(u_j)v^t_jp+pu_j\delta(v^t_j)p))
= \tau(\sum_j(v^t_jp\delta(u_j)+\delta(v^t_j)pu_j))\\
&=& \tau(\sum_j(v^t_j\delta(u_j)+\delta(v^t_j)u_j))
=\tau(\sum_j\delta(v^t_ju_j))=0.
\end{eqnarray*}
\end{proof}

Denote by $\tau_{\theta}$ the canonical trace on $A^{\infty}_{\theta}$ defined
by
\begin{eqnarray*} 
\tau_{\theta}(\sum_{h\in \Ze^n} c_hU_h)=c_0.
\end{eqnarray*}
Notice that, up to multiplication by a scalar, $\tau_{\theta}$ is the
unique continuous linear functional $\gamma$ on $A^{\infty}_{\theta}$ satisfying $\gamma\circ \delta_X=0$ for all
$X\in L^{\Ce}$. In the proof of the next lemma, which is trivial in the case $E$ a Hilbert bimodule,
we make crucial use of Theorem~\ref{cont of isom:thm}.

\begin{lemma} \label{trace:lemma}
Let $A^{\infty}_{\theta'}$, $E$, $A^{\infty}_{\theta}$, $\varphi$, and $\nabla$ be as in Theorem~\ref{complete me to orbit:thm}.
Denote by $\tau'$ the induced trace on $A^{\infty}_{\theta'}=\End(E_{A^{\infty}_{\theta}})$
obtained by the construction described above beginning with the trace $\tau_{\theta}$ on $A^{\infty}_{\theta}$.
Then $\tau'=\lambda \tau_{\theta'}$ for some $0\neq \lambda \in \Ce$.
\end{lemma}
\begin{proof}
By \cite[Corollary 2.3]{Schweitzer92} the subalgebra $M_k(A^{\infty}_{\theta})\subseteq M_k(A_{\theta})$ is closed under
the holomorphic functional calculus for any $k\in \Ne$. Therefore, by \cite[Proposition 4.6.2]{Blackadar},
every idempotent in $M_k(A^{\infty}_{\theta})$ is similar to a self-adjoint one, i.e.~a projection.
Consequently, in the definition of $\tau'$ we may choose
the idempotent $p\in M_k(A^{\infty}_{\theta})$ to be a projection.
By Example~\ref{cut down:example}, $pM_n(A^{\infty}_{\theta})p$ is closed under the holomorphic functional calculus
and has the Fr\'echet topology as the restriction of that on $M_k(A^{\infty}_{\theta})$.
By Theorem~\ref{cont of isom:thm}, applied to $A^{\infty}_{\theta'}$ and
the sub-$*$-algebra $(pM_k(A_{\theta})p)^{\infty}:=pM_k(A^{\infty}_{\theta})p$ of
$pM_k(A_{\theta})p$,
the natural isomorphism $A^{\infty}_{\theta'}=\End(E_{A^{\infty}_{\theta}})
\rightarrow pM_k(A^{\infty}_{\theta})p$ is a homeomorphism.
Therefore, the trace $\tau'$
is continuous on
$A^{\infty}_{\theta'}$ (as $\tau_{\theta}$ is obviously continuous on
$pM_k(A^{\infty}_{\theta})p$).
By Lemma~\ref{trace and der:lemma} we have
$\tau'\circ \delta_X=0$ for all $X\in L^{\Ce}$. Thus, (by the
remark above) $\tau'=\lambda \tau_{\theta'}$ for some $\lambda \in
\Ce$.

Using either that $p$ is full (since $E_{A^{\infty}_{\theta}}$ is a generator)
or that $\tau_{\theta}$ is faithful one sees that $\tau_{\theta}(p)>0$. Therefore
$\lambda\neq 0$.
\end{proof}

Let $A^{\infty}$ (resp. $B^{\infty}$) be a dense
sub-$*$-algebra of a C*-algebra $A$ (resp. $B$) closed under
the holomorphic functional calculus and equipped with a
Fr\'echet topology stronger than the C*-algebra norm topology.

\begin{lemma} \label{T:lemma}
The algebra $C^{\infty}(\Te, A^{\infty})$ is a dense sub-$*$-algebra of
the C*-algebra $C(\Te, A)$ closed under the holomorphic functional calculus and has a natural
Fr\'echet topology stronger than the C*-algebra norm topology.
If ${}_{B^{\infty}}E_{A^{\infty}}$ is a Morita equivalence bimodule, then
$C^{\infty}(\Te, B^{\infty})$ and $C^{\infty}(\Te, A^{\infty})$
are Morita equivalent with respect to the equivalence bimodule
${}_{C^{\infty}(\Te, B^{\infty})}C^{\infty}(\Te, E)_{C^{\infty}(\Te, A^{\infty})}$.
\end{lemma}
\begin{proof} Clearly $C^{\infty}(\Te, A^{\infty})$ is a sub-$*$-algebra of
$C(\Te, A)$. Since
it contains the algebraic tensor product $A^{\infty}\otimes C^{\infty}(\Te)$,
we see that $C^{\infty}(\Te, A^{\infty})$ is dense in $C(\Te, A)$.
Endow $C^{\infty}(\Te, A^{\infty})$ with the topology of uniform convergence on $\Te$ of the functions and of their
derivatives up to $s$ for every $s\in\Ne$. Clearly  this is a metrizable locally convex topology stronger than
the C*-algebra norm topology. We will show that this topology is complete.
Let $\{f_m\}_{m\in \Ne}$ be a Cauchy sequence in $C^{\infty}(\Te, A^{\infty})$.
Then the $s$-th derivatives $f^{(s)}_m$ converge uniformly to some continuous function $g_s:\Te\rightarrow A^{\infty}$.
Notice that $f^{(s)}_m(e^{iw})-f^{(s)}_m(e^{iv})=\int^w_vf^{(s+1)}_m(e^{it})\, dt$.
Taking limits we get $g_s(e^{iw})-g_s(e^{iv})=\int^w_vg_{s+1}(e^{it})\, dt$.
Consequently, $g'_s=g_{s+1}$. Thus $g_0\in C^{\infty}(\Te, A^{\infty})$.
It follows easily that $f_m\to g_0$ in $C^{\infty}(\Te, A^{\infty})$
as $m\to \infty$. So $C^{\infty}(\Te, A^{\infty})$ is complete.
Using the identity $a_1^{-1}-a_2^{-1}=a^{-1}_1(a_2-a_1)a^{-1}_2$ it is easy to see that
for any $f\in C^{\infty}(\Te, A^{\infty})$, if $f(t)$ is invertible in $A^{\infty}$ for
every $t\in \Te$, then $f^{-1}\in C^{\infty}(\Te, A^{\infty})$.
By \cite[Lemma 1.2]{Schweitzer92}, $f(t)$ is invertible in $A^{\infty}$ if and only
if it is invertible in $A$. Therefore, for any $f\in C^{\infty}(\Te, A^{\infty})$, if
it is invertible in $C(\Te, A)$ then it is invertible in $C^{\infty}(\Te, A^{\infty})$.
By \cite[Lemma 1.2]{Schweitzer92}, $C^{\infty}(\Te, A^{\infty})$ is closed under
the holomorphic functional calculus.

Let ${}_{B^{\infty}}E_{A^{\infty}}$ be a Morita equivalence bimodule.
Identify $E$ with $p(kA^{\infty})$ for some projection $p\in M_k(A^{\infty})$.
Define $P\in C^{\infty}(\Te, M_k(A^{\infty}))=M_k(C^{\infty}(\Te, A^{\infty}))$ to be the constant function
with value $p$ everywhere. Then
\begin{eqnarray*}
P^2=P,\quad  P(kC^{\infty}(\Te, A^{\infty}))=C^{\infty}(\Te, E),
\end{eqnarray*}
and
\begin{eqnarray*}
\End(C^{\infty}(\Te, E)_{C^{\infty}(\Te, A^{\infty})})&=&
PM_k(C^{\infty}(\Te, A^{\infty}))P \\
&=&C^{\infty}(\Te, pM_k(A^{\infty})p)=C^{\infty}(\Te, B^{\infty}).
\end{eqnarray*}
Since ${}_{B^{\infty}}E_{A^{\infty}}$ is a Morita equivalence
bimodule, the right module $E_{A^{\infty}}$ is a
generator,
which means that $A^{\infty}_{A^{\infty}}$ is a direct
summand of $rE_{A^{\infty}}$ for some $r\in \Ne$.
Equivalently,
there exist $\phi_j\in \Hom(E_{A^{\infty}}, A^{\infty})$ and $u_j\in E$
for $1\le j\le r$ such that $\sum_{1\le j\le
r}\phi_j(u_j)=1_{A^{\infty}}$. Denote by $\Phi_j:C(\Te, E)\rightarrow
C(\Te, A^{\infty})$ the map consisting of $\phi_j$ acting in fibres. Clearly,
$\Phi_j(C^{\infty}(\Te, E))\subseteq C^{\infty}(\Te, A^{\infty})$.
Let $\mathcal{U}_j\in C^{\infty}(\Te, E)$ denote the constant function
with value $u_j$ everywhere. Then $\sum_{1\le j\le
r}\Phi_j(\mathcal{U}_j)=1_{C^{\infty}(\Te,A^{\infty})}$.
Therefore
the right module $C^{\infty}(\Te, E)_{C^{\infty}(\Te, A^{\infty})}$ is a generator.
Similarly,
\begin{eqnarray*}
\End({}_{C^{\infty}(\Te, B^{\infty})}C^{\infty}(\Te,
E))=C^{\infty}(\Te, A^{\infty}),
\end{eqnarray*}
and ${}_{C^{\infty}(\Te, B^{\infty})}C^{\infty}(\Te,
E)$ is a finitely generated projective module and a generator.
Hence ${}_{C^{\infty}(\Te, B^{\infty})}C^{\infty}(\Te,
E)_{C^{\infty}(\Te, A^{\infty})}$ is a Morita equivalence bimodule.
\end{proof}

\begin{proposition} \label{K-groups:prop}
If ${}_{B^{\infty}}E_{A^{\infty}}$ is a Morita equivalence bimodule, then
there are natural group isomorphisms
$${\rm K}_0(B)\oplus {\rm K}_1(B)\rightarrow {\rm K}_0(C^{\infty}(\Te, B^{\infty}))
\rightarrow {\rm K}_0(C^{\infty}(\Te, A^{\infty})) \rightarrow {\rm K}_0(A)\oplus {\rm K}_1(A).
$$
\end{proposition}
\begin{proof}
By Bott periodicity we have a natural isomorphism
\begin{eqnarray*}
{\rm K}_0(A)\oplus {\rm K}_1(A)\rightarrow  {\rm K}_0(C(\Te, A)).
\end{eqnarray*}
By Lemma~\ref{T:lemma} the algebra $C^{\infty}(\Te, A^{\infty})$ is closed under the holomorphic functional calculus,
so we have
a natural isomorphism
$${\rm K}_0(C(\Te, A))\rightarrow {\rm K}_0(C^{\infty}(\Te, A^{\infty})).$$
Finally, the Morita equivalence bimodule
${}_{C^{\infty}(\Te, B^{\infty})}C^{\infty}(\Te, E)_{C^{\infty}(\Te, A^{\infty})}$ gives us
a natural
isomorphism $${\rm K}_0(C^{\infty}(\Te, B^{\infty}))
\rightarrow {\rm K}_0(C^{\infty}(\Te, A^{\infty})).$$
\end{proof}

We are ready to prove Theorem~\ref{complete me to orbit:thm}.

\begin{proof}[Proof of Theorem~\ref{complete me to orbit:thm}]
Consider the Fock space $\mathcal{F}^*=\Lambda((L^{\Ce})^*)$. Then
one can identify ${\rm K}_0(A_{\theta})$ and ${\rm K}_1(A_{\theta})$
with $\Lambda^{{\rm even}}(\Ze^n)$ and $\Lambda^{{\rm odd}}(\Ze^n)$
respectively \cite[Theorem 2.2]{Elliott84}. By
Proposition~\ref{K-groups:prop} we have a group isomorphism ${\rm
K}_0(A_{\theta'})\oplus {\rm K}_1(A_{\theta'})\rightarrow {\rm
K}_0(A_{\theta})\oplus {\rm K}_1(A_{\theta})$, which we shall think
of as $\psi:\Lambda(\Ze^n)\rightarrow \Lambda(\Ze^n)$. Of course
$\psi(\Lambda^{{\rm even}}(\Ze^n))=\Lambda^{{\rm even}}(\Ze^n)$.
Notice that $L^*$ acts on $\mathcal{F}^*$ via multiplication. Also
$L$ acts on $\mathcal{F}^*$ via contraction. Let $a^1, {\cdots},
a^n$ and $b_1, {\cdots}, b_n$ denote the standard bases of $L^*$ and
$L$ respectively. Denote $(a^1, {\cdots}, a^n)$ and $(b_1, {\cdots},
b_n)$ by $\vec{a}$ and $\vec{b}$ respectively. Denote by
$\mathscr{A}$ the matrix of $\varphi$ with respect to $\vec{b}$.
Denote by $\Phi$ the $n\times n$ matrix $\frac{1}{2\pi
i}([\nabla_{b_j}, \nabla_{b_k}])$. Using the Chern character which
was defined in \cite{Connes80} and calculated for noncommutative
tori in \cite{Elliott84}, Schwarz showed that in the case of
complete Morita equivalence (not assumed here)
 the matrix
\begin{eqnarray} \label{g:eq}
\quad \quad \quad g
&=&\begin{pmatrix} S & R \\ N &  M \end{pmatrix} \\
\nonumber &:=&\begin{pmatrix}
\mathscr{A}^{-1}+\theta\Phi\mathscr{A}^{-1} &
-\mathscr{A}^{-1}\theta'-\theta\Phi\mathscr{A}^{-1}\theta'+\theta\mathscr{A}^t
\\ \Phi\mathscr{A}^{-1} &
-\Phi\mathscr{A}^{-1}\theta'+\mathscr{A}^t
\end{pmatrix}
\end{eqnarray}
is in $\SO(n, n|\Re)$ and there is a linear operator $V$ on
$\mathcal{F}^*$ extending $\psi|_{\Lambda^{{\rm even}}(\Ze^n)}$ such
that
\begin{eqnarray} \label{W:eq}
V(\vec{b}, \vec{a})V^{-1}=(\vec{b}, \vec{a})g.
\end{eqnarray}
Our equations (\ref{g:eq}) and (\ref{W:eq}) are exactly the
equations (49), (50) and (53) of \cite{Schwarz98}, in slightly
different form. From the equation (\ref{g:eq})
above, Schwarz deduced
\begin{eqnarray} \label{action2:eq}
\theta=(S\theta'+R)(N\theta'+M)^{-1},
\end{eqnarray}
which is our desired conclusion, except for the assertion that the
matrix $g$ belongs to $M_{2n}(\Ze)$ (and hence to $\SO(n, n|\,
\Ze)$). Note that although in the definition of the Chern character
in \cite{Connes80} Connes required the connections to be Hermitian,
all the arguments there hold for arbitrary connections. Using
Lemma~\ref{trace:lemma} one checks that Schwarz's argument to get
(\ref{W:eq}) and (\ref{action2:eq}) still works in our situation (in
which neither the connection nor the Morita equivalence are
Hermitian) except that now we can only say that $V$ acts on
$\mathcal{F}^*$ and $g$ is in $\SO(n, n|\, \Ce)$; in other words,
$g$ might not be in $M_{2n}(\Re)$ a priori. In the complete Morita
equivalence case, referring to the fact that $V$ maps $\Lambda^{{\rm
even}}(\Ze^n)$ into itself and satisfies (\ref{W:eq}) with $g\in
\SO(n, n|\, \Re)$, Schwarz concluded that $g\in M_{2n}(\Ze)$ for the
case $n>2$ (this is not true for $n=2$), so that $g\in \SO(n, n|\,
\Ze)$, as desired. We have not been able to understand this part of
the argument, and so we shall follow another route: we assert that
actually $V$ extends all of $\psi$ and hence maps all of
$\Lambda(\Ze^n)$ onto itself (not just $\Lambda^{{\rm
even}}(\Ze^n)$).
 This follows from applying
Schwarz's argument to the Morita equivalence bimodule
${}_{C^{\infty}(\Te, A^{\infty}_{\theta'})}C^{\infty}(\Te,
E)_{C^{\infty}(\Te, A^{\infty}_{\theta})}$ (combining even and odd
degrees) of Lemma~\ref{T:lemma} instead of to
${}_{A^{\infty}_{\theta'}}E_{A^{\infty}_{\theta}}$. For the
convenience of the reader, we sketch the argument here.

We recall the definition of the Chern character first. Consider the
trivial Lie algebra $L^{\Ce}\oplus \Ce$. It acts on $C^{\infty}(\Te,
A^{\infty}_{\theta'})$ as derivations $\delta$ by extending the
action of $L^{\Ce}$ on $A^{\infty}_{\theta'}$, with $L^{\Ce}$ acting
on $C^{\infty}(\Te)$ trivially,  and the action of the canonical
unit vector of $\Ce$ being the differentiation with respect to the
anti-clockwise unit vector field on $\Te$. Tensoring
$\tau_{\theta'}$ with the Lebesgue integral on $C^{\infty}(\Te)$, we
obtain an $L^{\Ce}\oplus \Ce$-invariant trace on $C^{\infty}(\Te,
A^{\infty}_{\theta'})$, which we still denote by $\tau_{\theta'}$.
For a finitely generated projective right $C^{\infty}(\Te,
A^{\infty}_{\theta'})$-module $F'_{C^{\infty}(\Te,
A^{\infty}_{\theta'})}$, the Chern character $\ch F'$ is defined by
\begin{eqnarray} \label{Chern:eq}
\ch F'&=&\tau_{\theta'}(e^{\Omega'/2\pi
i})=\sum_{j=0}\frac{1}{j!}\tau_{\theta'}((\Omega')^j)\cdot
\frac{1}{(2\pi i)^j}\\
\notag &\in& \Lambda^{\rm even}((L^{\Ce}\oplus
\Ce)^*)=\Lambda((L^{\Ce})^*),
\end{eqnarray}
where $\Omega'\in \End(F'_{C^{\infty}(\Te,
A^{\infty}_{\theta'})})\otimes \Lambda^2((L^{\Ce}\oplus \Ce)^*)$ is
the curvature of an arbitrary connection on $F'$ (with respect to
the action of $L^{\Ce}\oplus \Ce$ on $C^{\infty}(\Te,
A^{\infty}_{\theta'}))$ and we have extended $\tau_{\theta'}$ to
$\End(F'_{C^{\infty}(\Te, A^{\infty}_{\theta'})})$ as in the
paragraph before Lemma~\ref{trace and der:lemma}. This determines
the Chern character $\ch:\Lambda(\Ze^n)={\rm K}_0(A_{\theta'})\oplus
{\rm K}_1(A_{\theta'})={\rm K}_0(C^{\infty}(\Te,
A^{\infty}_{\theta'}))\rightarrow \Lambda((L^{\Ce})^*)$ as a group
homomorphism. The Chern character ${\rm K}_0(C^{\infty}(\Te,
A^{\infty}_{\theta}))\rightarrow \Lambda((L^{\Ce})^*)$ is defined
similarly.

Now let us consider the finitely generated right $C^{\infty}(\Te,
A^{\infty}_{\theta})$-module $F:=F'\otimes_{C^{\infty}(\Te,
A^{\infty}_{\theta'})}E$. To get a connection of $F_{C^{\infty}(\Te,
A^{\infty}_{\theta})}$ from that of $F'_{C^{\infty}(\Te,
A^{\infty}_{\theta'})}$, let us extend $\varphi$ to $L^{\Ce}\oplus
\Ce\rightarrow L^{\Ce}\oplus \Ce$ as simply being the identity map
on $\Ce$, and also extend $\nabla$ to $L^{\Ce}\oplus \Ce\rightarrow
\Hom_{\Ce}(C^{\infty}(\Te, E))$ in such a way that the action of
$L^{\Ce}$ of $C^{\infty}(\Te, E)$ is fibrewise the original $\nabla$
and the action of the canonical unit vector of $\Ce$ on
$C^{\infty}(\Te, E)$ is the differentiation with respect to the
anti-clockwise unit vector field on $\Te$. Then $\nabla_X$  is a
connection for both $(C^{\infty}(\Te, E)_{C^{\infty}(\Te,
A^{\infty}_{\theta})}, \delta_X)$ and $({}_{C^{\infty}(\Te,
A^{\infty}_{\theta'})}C^{\infty}(\Te, E), \delta_{\varphi(X)})$ for
every $X\in L^{\Ce}\oplus \Ce$, and furthermore $\nabla$ has
constant curvature $\pi i\sum_{j, k}\Phi_{jk}a^j\wedge a^k$ in
$\End(C^{\infty}(\Te, E)_{C^{\infty}(\Te,
A^{\infty}_{\theta})})\otimes \Lambda^2((L^{\Ce}\oplus \Ce)^*)$.

For any connection $\nabla'_{\varphi(X)}$ of $(F'_{C^{\infty}(\Te,
A^{\infty}_{\theta'})}, \delta_{\varphi(X)})$, one checks readily
that $\nabla'_{\varphi(X)}\otimes \id+\id\otimes \nabla_X$ is a
connection of $(F_{C^{\infty}(\Te, A^{\infty}_{\theta})},
\delta_{X})$. If we choose the connections of $(F_{C^{\infty}(\Te,
A^{\infty}_{\theta})}, \delta)$ in this way, then the curvature
$\Omega$ is calculated by
\begin{eqnarray} \label{curvature:eq}
\Omega=\varphi^*(\Omega')+\pi i\sum_{j, k}\Phi_{jk}a^j\wedge a^k,
\end{eqnarray}
where $\varphi^*$ denotes the linear isomorphism
$\Lambda((L^{\Ce}\oplus \Ce)^*)\rightarrow \Lambda((L^{\Ce}\oplus
\Ce)^*)$ induced by $\varphi$ and we identify
$\End(F'_{C^{\infty}(\Te, A^{\infty}_{\theta'})})$ with a subalgebra
of $\End(F_{C^{\infty}(\Te, A^{\infty}_{\theta})})$ via identifying
$T\in \End(F'_{C^{\infty}(\Te, A^{\infty}_{\theta'})})$ with
$T\otimes \id$. Let $\lambda$ be the constant in
Lemma~\ref{trace:lemma}. Since $\tau_{\theta'}$ (resp.
$\tau_{\theta}$) was extended to $C^{\infty}(\Te,
A^{\infty}_{\theta'})$ (resp. $C^{\infty}(\Te,
A^{\infty}_{\theta})$) via tensoring with the Lebesgue integral on
$\Te$, the conclusion of Lemma~\ref{trace:lemma} actually holds on
$C^{\infty}(\Te, A^{\infty}_{\theta'})=\End(C^{\infty}(\Te,
E)_{C^{\infty}(\Te, A^{\infty}_{\theta})})$. Thus
$\tau_{\theta}(a)=\lambda \tau_{\theta'}(a)$ for all $a\in
\End(F'_{C^{\infty}(\Te, A^{\infty}_{\theta'})})$. Therefore, we
have
\begin{eqnarray} \label{curvature 2:eq}
\ch F&=&\tau(e^{\Omega/2\pi i})=\tau(e^{\frac{1}{2}\sum_{j,
k}\Phi_{jk}a^j\wedge a^k}\varphi(e^{\Omega'/2\pi i}))\\
\notag &=& \lambda e^{\frac{1}{2}\sum_{j,
k}\Phi_{jk}a^ja^k}\varphi^*(\ch F').
\end{eqnarray}

Denote by $\mu(F')$ the equivalence class of $F'$ in ${\rm
K}_0(C^{\infty}(\Te, A^{\infty}_{\theta'}))=\Lambda(\Ze^n)$. By
Theorem 4.2 of \cite{Elliott84} (the sign there must be reversed;
see page 137 of \cite{DEKR85}), we have
\begin{eqnarray} \label{Chern 2:eq}
\ch F'=e^{-\frac{1}{2}\sum_{j, k}\theta'_{jk}b_jb_k}\mu(F').
\end{eqnarray}
Combining the equations (\ref{Chern:eq}), (\ref{curvature 2:eq}),
and (\ref{Chern 2:eq}) together, we see that the map $\psi$ is the
restriction of the linear operator $V:=V_1V_2V_3V_4\in
\Hom_{\Ce}(\mathcal{F}^*)$ on $\Lambda(\Ze^n)$, where
\begin{eqnarray*}
V_1f &=& e^{\frac{1}{2}\sum_{j, k}\theta_{jk}b_jb_k}f,
\\
V_2f &=& \lambda e^{\frac{1}{2}\sum_{j,
k}\Phi_{jk}a^ja^k}
f, \\
V_3f &=& \varphi^*(f), \\
V_4f &=& e^{-\frac{1}{2}\sum_{j, k}\theta'_{jk}b_jb_k}f,
\end{eqnarray*}
for $f\in \mathcal{F}^*$.

 Each linear operator $V_k$ is a {\it linear
canonical transformation} in the sense that
$$ V_k(\vec{b}, \vec{a})V^{-1}_k=(\vec{b}, \vec{a})g_k $$
for some $g_k\in M_{2n}(\Ce)$. In fact, a simple calculation yields
\begin{eqnarray*}
g_1=\begin{pmatrix} I & \theta \\ 0 &  I \end{pmatrix}, \quad
g_2=\begin{pmatrix} I & 0 \\ \Phi &  I \end{pmatrix}, \quad
g_3=\begin{pmatrix} \mathscr{A}^{-1} & 0 \\ 0 &  \mathscr{A}^t
\end{pmatrix}, \quad g_4=\begin{pmatrix} I & -\theta' \\ 0 &  I
\end{pmatrix}.
\end{eqnarray*}
Set $g=g_1g_2g_3g_4$. Then the equations (\ref{g:eq}) and
(\ref{W:eq}) hold.

Notice that each $g_k$ belongs to $\SO(n, n|\, \Ce)$, i.e. it
satisfies the equations (\ref{O(n,n|R):eq}) and has determinant $1$.
Hence so also does $g$.

Since $V$ extends the automorphism $\psi$ of $\Lambda(\Ze^n)$, $g$
is easily seen to be in $M_{2n}(\Ze)$ by applying (\ref{W:eq}) to
the canonical vectors $1$ and $a^j, \, 1\le j\le n$, in the Fock
space $\mathcal{F}^*$. Therefore, $g$ belongs to $\SO(n, n|\, \Ze)$.
\end{proof}

\begin{remark} \label{isomorphism:remark}
Let us indicate briefly how the proof of the part (2) of
Theorem~\ref{Morita equiv:thm} leads to a new proof of the main
result of \cite{CEGJ85}, namely, if $\theta'$ and $\theta$ are in
$\mathcal{T}'_n\cap \cT^{\flat}_n$ and the algebras
$A^{\infty}_{\theta'}$ and $A^{\infty}_{\theta}$ are isomorphic,
then the bicharacters $\rho_{\theta'}$ and $\rho_{\theta}$ of
$\Ze^n$ are isomorphic. On using the given isomorphism
$A^{\infty}_{\theta'}\rightarrow A^{\infty}_{\theta}$, the vector
space $E=A^{\infty}_{\theta}$ becomes a Morita equivalence bimodule
for $A^{\infty}_{\theta'}$ and $A^{\infty}_{\theta}$ in a natural
way. The Chern character \cite{Connes80} of the free module
$E_{A^{\infty}_{\theta}}$($={A^{\infty}_{\theta}}_{A^{\infty}_{\theta}}$)
is $1$. Therefore, the constant curvature connection of
Theorem~\ref{me to complete:thm} has in fact curvature zero. It
follows that we have $N=0$ in (\ref{g:eq}). Since $g\in \SO(n, n|\,
\Ze)$, the block entry $S$ must belong to $\GL(n, \Ze)$. A simple
calculation (which is trivial in the case also $R=0$) shows that the
bicharacters associated to $\theta'$ and $\theta=g\theta'$ are
isomorphic (by means of the automorphism $S$ of $\Ze^n$).
\end{remark}

\end{document}